\documentclass[11pt,a4paper]{article}
\usepackage[utf8]{inputenc}
\usepackage{amsmath}
\usepackage{amsfonts}
\usepackage{amssymb}
\usepackage{makeidx}
\usepackage{graphicx}
\usepackage{mathabx} 
\usepackage{dsfont} 
\usepackage[left=2cm,right=2cm,top=2cm,bottom=2cm]{geometry}
\usepackage{color}
\usepackage{verbatim} 
\usepackage{tocbibind} 
\usepackage{hyperref} 
\usepackage{mathtools} 

\newtheorem{satz}{Lemma}[section]
\newtheorem{theorem}{Theorem}[section]
\newtheorem{definition}{Definition}[section]

\newtheorem{proposition}{Proposition}[section]
\newtheorem{corollary}{Corollary}[section]

\begin{document}

\setlength{\parindent}{0pt} 


\ \vspace{3cm}

\begin{center}
\Large\textbf{FLUID-RIGID BODY INTERACTION IN A COMPRESSIBLE ELECTRICALLY CONDUCTING FLUID} \\[10mm] \large{JAN SCHERZ}$^{1,2,3}$ \\[10mm]
\end{center}

\begin{itemize}
\item[$^1$] Department of Mathematical Analysis, Faculty of Mathematics and Physics, Charles University in Prague, Sokolovská 83, Prague 8, 18675, Czech Republic
\item[$^2$] Mathematical Institute, Academy of Sciences, \v{Z}itná 25, Prague 1, 11567, Czech Republic
\item[$^3$] Institute of Mathematics, University of Würzburg, Emil-Fischer-Str. 40, 97074 Würzburg, Germany
\end{itemize}   

\bigskip

\begin{center}
\Large\textbf{Abstract} \\[4mm]
\end{center}
We consider a system of multiple insulating rigid bodies moving inside of an electrically conducting compressible fluid. In this system we take into account the interaction of the fluid with the bodies as well as with the electromagnetic fields trespassing both the fluid and the solids. The main result of this article yields the existence of weak solutions to the system. While the mechanical part of the problem can be dealt with via a classical penalization method, the electromagnetic part requires an approximation by means of a hybrid discrete-continuous in time system: The discrete part of the approximation enables us to handle the solution-dependent test functions in our variational formulation of the induction equation, whereas the continuous part makes sure that the non-negativity of the density and subsequently a meaningful energy inequality is preserved in the approximate system.

{\centering \section{Introduction} \label{Introduction} \par }

The goal of this article is the proof of the existence of weak solutions to a system of PDEs modelling the motion of multiple non-conducting rigid objects inside of an electrically conducting compressible fluid taking the interplay with the electromagnetic fields in these materials into account. The problem can be regarded as belonging to both the research areas of \textit{magnetohydrodynamics (MHD)} and \textit{fluid-structure interaction (FSI)}. Indeed, MHD (c.f.\ \cite{cabannes,davidson,kl}), on the one hand, describes a coupling between the Maxwell system and the Navier-Stokes equations which models the influence of the electromagnetic fields on the electrically conducting fluid and vice versa. FSI, on the other hand, models the interaction between the fluid and the solid bodies through a coupling between the Navier-Stokes system and the balance of mass and momentum of the rigid bodies. 

The interplay between the electromagnetic fields and the solids occurs implicitly via the fluid due to the non-conductivity of the solid material. This paper represents the extension of \cite{incompressiblecase}, where the corresponding problem was studied and solved for incompressible fluids, to the compressible case. Potential applications can be found in the area of biomechanics. A specific example is constituted by capsule endoscopy, a medical procedure during which small cameradevices are sent through the body with the aim of detecting diseases. Such endocapsules of a microscopic scale can be navigated through the (electrically conducting) blood by controlling their movement via the application of electromagnetic forces, c.f.\ \cite[Section 4.4]{drugdelivery}, \cite{endocapsules}. This technique can also be applied in the problem of drug delivery, in which microrobots are constructed to deliver drugs directly to the targeted area of the body without affecting healthy tissue. However, in view of the usage of electrically conducting microrobots for these purposes an extension of the current model to electrically conducting rigid bodies would be needed.

In order to classify our result within the wide range of works in MHD and FSI, we give a brief summary of the related literature. The existence of weak solutions to the MHD problem for a compressible fluid but without any rigid bodies involved is for example shown in \cite{sart} and \cite{blancducomet}. The fluids considered in the latter one of these articles are electrically as well as thermally conducting and moreover the existence of strong solutions is addressed therein. The existence of weak solutions to the MHD problem in the incompressible case is proved in \cite{gbl}.

On the FSI side of the problem we first mention \cite{G2,SER3} for a general introduction to the fluid-rigid body interaction problem. Early existence results for weak solutions to this problem were obtained mainly in the incompressible case and include the articles \cite{CST,desjardinsesteban,GLSE,HOST}, wherein the proof of local in time existence, i.e. existence up to a contact between several bodies or a body and the domain boundary, is achieved in the case of two and three spatial dimensions. The corresponding problem in the compressible case was considered in \cite{compressibledesjardinsesteban}. A proof of the global in time existence of weak solutions to a model describing the interaction between multiple rigid bodies and a compressible fluid is given in \cite{feireisl}. Corresponding results are also available for incompressible fluids, c.f.\ \cite{tucsnak} for the two-dimensional and \cite{incompressiblefeireisl} for the three-dimensional case. The article \cite{tucsnak} in particular touches upon the question whether contacts between the bodies with each other or the domain boundary are possible and shows that such collisions can only occur if the relative acceleration and velocity between the colliding objects vanish. Moreover, the problem has been studied in the case of the Navier-slip instead of the classical no-slip boundary condition, c.f.\ \cite{compressiblecase}, and also the question about the existence of strong solutions has been investigated in both the compressible, c.f.\ \cite{boulakiaguerrero,haakmaitytakahashi,hiebermurata,roytakahashi}, and the incompressible case, c.f.\ \cite{GGH13,T,Wa}.

The articles \cite{guermondminev2d} (for the case of two spatial dimensions) and \cite{guermondminev} (for the $3$D case) can be considered as a first step towards the coupling between the MHD and the FSI problem. The authors thereof studied a model of an incompressible electrically conducting fluid flowing around a fixed non-conducting solid region. This model served as the basis for the electromagnetic part of the model in \cite{incompressiblecase}, wherein Benešová, Nečasová, Schlömerkemper and the author of this article showed the (local in time) existence of weak solutions to the problem of one movable insulating rigid object travelling through an electrically conducting incompressible fluid. The main goal of the present article is to extend the latter result to the case of a compressible fluid. More precisely, we are able to prove the global in time existence of weak solutions to the interaction problem of finitely many insulating rigid bodies, an electrically conducting compressible fluid and the electromagnetic fields trespassing these materials, see Section \ref{mainresultsection}.

As in the incompressible case in \cite{incompressiblecase}, the main difficulty in the proof of the existence result is caused by the test functions in the weak formulation of the induction equation, c.f.\ \eqref{-340}, \eqref{-218} below, which are chosen such that they depend on the solid region. While such test functions do not generate any problems in the case of an immovable solid region (see e.g.\ the proofs of \cite[Theorem 2.1]{guermondminev2d} and \cite[Theorem 2.3]{guermondminev}), difficulties arise in our scenario, where the solid domain depends on the overall solution to the system, which causes our problem to be highly coupled. Following the proof in the incompressible situation, we thus make use of a time discretization, which allows us to deal with this problem by decoupling it: At each discrete time we first calculate the domain of the solid bodies, which suggests a suitable definition for the test functions in the induction equation at that specific time. Only after this we solve the induction equation itself, which can then be achieved via standard methods. In the compressible situation, however, this procedure turns out to cause more problems. In particular, the author could not find a suitable way to discretize the Navier-Stokes system while preserving the non-negativity of the density. This is essential to obtain the uniform bounds from the energy inequality required to pass to the limit in the approximate system. We handle this problem by choosing a hybrid approximation system, in which the induction equation is discretized in time via the Rothe method (\cite[Section 8]{roubicek}), whereas the mechanical part of the system is studied as a time-dependent problem on the small intervals between the discrete time points. The non-negativity of the density can then be derived by classical arguments and, by choosing the coupling terms in a suitable way, the discrete electromagnetic part and the continuous mechanical part of the system can be combined to an energy inequality with all desired features. To a smaller extent, we already used such a hybrid approximation in the incompressible situation \cite{incompressiblecase}, where, however, the time dependency was restricted to the transport equation for the characteristic function of the solid region. The expansion of this idea to the whole mechanical part of the system, in order to deal with the problems outlined above, is what constitutes the main novelty in the proof of our main result.

The problem of the solution dependent test functions also appears in the weak formulation of the momentum equation, c.f.\ \eqref{-339}, \eqref{-217} below. In this situation, however, we have the penalization method used e.g.\ in \cite{feireisl} and \cite{tucsnak} readily available which allows us to evade the problem. More precisely, in this penalization method a sequence of approximate solutions to some fluid-only problems with classical test functions is constructed. Passing to the limit in this approximation one then returns to a fluid-rigid body interaction system by letting the viscosity of the fluid rise to infinity in the later solid region.

The paper is structured as follows: We begin by summarizing the notation needed for our model in Section \ref{model} and subsequently present the model, divided into a mechanical and an electromagnetic part, in Sections \ref{mechanicalsubsystem} and \ref{electromagneticsubsystem}. After introducing some additional notation in Section \ref{notation}, we present the variational formulation of the above model in Section \ref{weaksolutionssection} as well as our main result in Section \ref{mainresultsection}. In Section \ref{approximatesystem} we explain the main ideas for the proof of this result, which is based on an approximation of the original system. In Section \ref{existencediscretesolution} we solve this approximate problem and finally, in the remaining Sections \ref{deltatlimit}--\ref{alphalimit}, we pass to the limit in the approximation, proving the existence of a weak solution to the original system.\\

\subsection{Model} \label{model}

The model we consider describes several non-conducting rigid bodies travelling through an electrically conducting compressible fluid as well as the involved electromagnetic fields. This model is a combination of (i) the mechanical fluid-rigid body interaction model used in \cite{feireisl} and (ii) the Maxwell system in the model used in \cite{guermondminev2d,guermondminev} for the description of an electrically conducting fluid surrounding an immovable, non-conducting solid region. It is further the extension of the corresponding model for the incompressible case in \cite{incompressiblecase} to the compressible situation. Let $T>0$ and let $\Omega \subset \mathbb{R}^3$ be a bounded domain. Inside of $\Omega$ we consider $m \in \mathbb{N}$ insulating rigid bodies, the positions of which at time $t \in [0,T]$ are described through subsets $S^i(t) \subset \Omega$, $i=1,...,m$. The complement of the solid domain,
\begin{align}
F(t):= \Omega \setminus \overline{S}(t),\quad \quad \text{with}\quad S(t) := \bigcup_{i=1}^m S^i(t), \nonumber
\end{align}

contains an electrically conducting viscous nonhomogeneous compressible fluid. We denote by $Q$ the time-space domain $Q:=(0,T)\times \Omega$, which we split into a solid part $Q^{s}$ and a fluid part $Q^f$,
\begin{align}
Q^{s}=Q^{s}\left(S\right) := \left\lbrace (t,x) \in Q:\ x \in S(t) \right\rbrace,\quad \quad Q^f=Q^f(S) := \left\lbrace (t,x) \in Q:\ x \in F(t) \right\rbrace. \label{-353}
\end{align}

For any function defined on $Q$ we mark its restriction to $Q^f$ or $Q^s$ by the superscript $f$ or $s$, respectively.

\begin{figure}[h!]
\begin{center}
\includegraphics[scale=0.5]{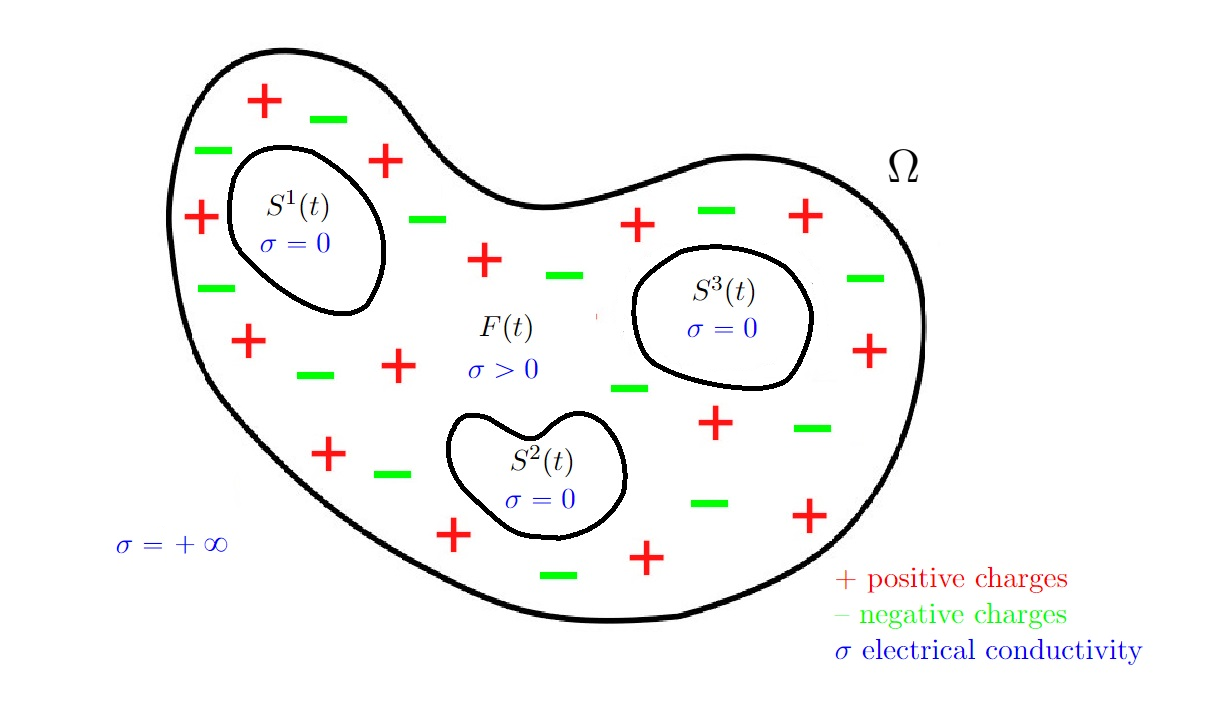} 
\caption{A domain filled with an electrically conducting compressible fluid and three insulating rigid bodies.} \label{figure}
\end{center}
\end{figure}

The interaction between the fluid, the solids and the electromagnetic fields in the domain is characterized through the mass density $\rho:Q \rightarrow \mathbb{R}$, the velocity field $u:Q \rightarrow \mathbb{R}^3$, the magnetic induction $B:Q \rightarrow \mathbb{R}^3$, the electric field $E:Q \rightarrow \mathbb{R}^3$ and the electric current $j:Q \rightarrow \mathbb{R}^3$. As indicated above, our overall model, which determines these functions, can be divided into two subsystems: The mechanical subsystem for the description of $\rho$ and $u$ and the electromagnetic subsystem for the description of $B$, $E$ and $j$.

\subsection{Mechanical subsystem} \label{mechanicalsubsystem}

The mechanical quantities evolve according to the compressible Navier-Stokes equations in the fluid domain and the balance of linear and angular momentum of the rigid bodies in the solid region, respectively,
\begin{align}
\partial _t \rho + \operatorname{div} (\rho u) = 0 \quad \quad &\text{in } Q^f, \label{526} \\
\partial _{t} (\rho u) + \text{div} (\rho u \otimes u) + \nabla p = \text{div} \mathbb{T} + \rho g + \frac{1}{\mu} \text{curl}B\times B\quad \quad &\text{in } Q^f, \label{527} \\
m^i \frac{d}{dt} V^i(t) = \frac{d}{dt} \int_{S^i(t)} \rho u\ dx = \int_{\partial S^i(t)} \left[\mathbb{T} - p\text{Id}\right] \cdot \text{n}\ d\sigma + \int _{S^i(t)} \rho g \ dx,& \quad t \in [0,T],\ i = 1,...,m, \label{528} \\
\frac{d}{dt}\left( \mathbb{J}^i(t)w^i(t) \right) = \frac{d}{dt} \int _{S^i(t)} \rho \left(x - h^i\right) \times u \ dx\ & \nonumber \\
= \int_{\partial S^i(t)} (x - h^i) \times \left[\mathbb{T} - p\text{Id} \right] \text{n}\ d\sigma + \int _{S^i(t)} \rho \left(x - h^i\right) \times g\ dx,& \quad t \in [0,T],\ i=1,...,m, \label{529}
\end{align}

combined with the boundary and interface conditions
\begin{align}
u(t) &= 0\quad \text{on } \partial \Omega,\quad \quad &&\ \ \ \ \ \ \ \ \ \ u^f(t) - u^{s}(t) = 0 \quad \text{on } \partial S(t). \label{833}
\end{align}

The identity \eqref{526} is known as the continuity equation. In the momentum equation \eqref{527} we see the pressure $p$, for which we assume the isentropic constitutive relation
\begin{align}
p = p \left(\rho^f \right) = a \left( \rho^f \right)^\gamma, \quad a>0,\ \gamma > \frac{3}{2}, \nonumber
\end{align}

and the stress tensor
\begin{align}
\mathbb{T} = \mathbb{T}(u) := 2\nu \mathbb{D}(u) + \lambda \text{Id}\ \operatorname{div} u,\quad \quad \mathbb{D}(u) := \frac{1}{2} \nabla u + \frac{1}{2} (\nabla u)^T \nonumber
\end{align}

with the viscosity coefficients $\nu,\lambda \in \mathbb{R}$ which satisfy $\nu > 0,\ \lambda + \nu \geq 0$. Moreover, we have two forcing terms: The external force $g$ and the Lorentz force $\frac{1}{\mu} \text{curl}B\times B$ with the magnetic permeability $\mu> 0$. We assume that
\begin{align}
\mu \ \text{is constant in the whole domain} \ Q. \nonumber
\end{align}

Since, in general, the magnetic permeability takes different values in conducting and insulating materials, this assumption is physically restrictive. However, it is necessary for the transition conditions on the magnetic induction $B$, c.f.\ Section \ref{electromagneticsubsystem} below. The Lorentz force constitutes the connection of the mechanical to the electromagnetic subsystem. In the balance of mass \eqref{528} and momentum \eqref{529} of the rigid bodies instead the Lorentz force does not appear, which is in accordance with the assumption of the bodies being non-conducting. The notation used in these relations includes the total mass $m^i$ of the $i$-th body, its center of mass $h^i$ and the associated inertia tensor $\mathbb{J}^i$,
\begin{align}
m^i := \int _{S^i(t)} \rho (t,x)\ dx,\quad h^i(t):= \frac{1}{m^i}\int _{S^i(t)} \rho (t,x)x\ dx, \quad \quad &t \in [0,T], \nonumber \\
\mathbb{J}^i(t)a \cdot b := \int _{S^i(t)}\rho (t,x)\left[ a \times \left( x - h^i(t) \right) \right] \cdot \left[ b \times \left( x - h^i(t) \right) \right]\ dx,\quad \quad &t \in [0,T],\quad a,b \in \mathbb{R}^3. \nonumber
\end{align}

The equations \eqref{528} and \eqref{529} determine the translational and rotational velocities $V^i$ and $w^i$ of the $i$-th rigid body, respectively, allowing us to express its total velocity as a rigid velocity field
\begin{align}
u(t,x) = u^{s^i}(t,x) := V^i(t) + w^i(t) \times \left( x - h^i(t) \right) \quad \quad \text{for } t \in [0,T],\ x \in S^i(t). \nonumber
\end{align}

The coupling between the fluid and the solids is incorporated into the surface integrals in (\ref{528}) and (\ref{529}) and into the no-slip interface condition in \eqref{833}. Indeed, on the one hand, the presence of the Cauchy-stress $\mathbb{T} - p\text{Id}$ of the fluid in (\ref{528}) and (\ref{529}) shows how the velocity and the pressure of the fluid affect the motion of the bodies. On the other hand, considering the Navier-Stokes system \eqref{526}, \eqref{527}, \eqref{833} by itself, one can regard the interface condition in \eqref{833} as a part of the boundary conditions, which describes the impact of the velocity of each solid on the fluid velocity. The first identity in \eqref{833} is the standard no-slip boundary condition.

\subsection{Electromagnetic subsystem} \label{electromagneticsubsystem}

The electromagnetic quantities are determined by the following version of the Maxwell system,
\begin{align} \frac{1}{\mu} \text{curl} B &= \left\{
                \begin{array}{ll}
                  j + J \ \ \ &\text{in } Q^f\\
                  0 &\text{in } Q^s
                \end{array}
              \right., \label{94} \\
\partial_t B + \text{curl} E &= 0 \quad \quad \quad \quad \quad \text{in } Q^f \ \text{and } Q^s, \label{95} \\
\text{div} E &= 0 \quad \quad \quad \quad \quad \text{in } Q^s, \label{129} \\
\text{div} B &= 0 \quad \quad \quad \quad \quad \text{in } Q^f\ \text{and } Q^s,\label{876}
\end{align}

together with Ohm's law
\begin{align}
&j = \sigma (E + u \times B) \quad \text{in } Q,\quad \sigma = \left\{
                \begin{array}{ll}
                  \sigma ^f > 0 \ \ \ &\text{in } Q^f\\
                  \sigma ^s = 0 &\text{in } Q^s
                \end{array}
              \right. \label{877}
\end{align}

and the boundary and interface conditions
\begin{align}
B(t) \cdot \text{n} &= 0\quad \text{on } \partial \Omega,\quad \quad &&\ \ \ \ \ \ \ \ \ B^f(t) - B^{s}(t) = 0 \quad \text{on } \partial S(t), \label{831} \\
E(t) \times \text{n} &= 0 \quad \text{on } \partial \Omega,\quad \quad &&\left(E^f(t) - E^{s}(t)\right) \times \text{n} = 0 \quad \text{on } \partial S(t). \label{832}
\end{align}

In this system we have Ampère's law \eqref{94}, the Maxwell-Faraday equation (\ref{95}), Gauss's law (\ref{129}) and Gauss's law for magnetism (\ref{876}). In comparison to the general situation, these equations have undergone two kinds of reductions, c.f.\ \cite{guermondminev2d}, \cite{guermondminev}: First, in the solid region, the system is adjusted to the assumption of the rigid bodies being insulating and second, in the fluid domain, the system is reduced according to the \textit{magnetohydrodynamic approximation}, see \cite{davidson, eringenmaugin}. A justification for the latter simplification, which is inherent to magnetohydrodynamics, from a physical point of view is for example given in \cite{jiliang,jiliang2}. The mathematical use of the magnetohydrodynamic approximation consists of the fact that it enables us to summarize the whole electromagnetic problem in $Q^f$ into a problem for the magnetic induction $B$, c.f.\ the induction equation \eqref{-218} in our definition of variational solutions below. Once $B$ has been determined, the remaining unknowns $E$ and $j$ are given explicitly through the relations \eqref{94} and \eqref{877}. The quantity $J$ in \eqref{94} represents, as in \cite{guermondminev2d}, \cite{guermondminev}, a supplementary external force. Ohm's law \eqref{877} is what links the electromagnetic subsystem to the mechanical subsystem \eqref{526}--\eqref{833} by describing the influence of the (fluid) velocity $u=u^f$ on the electromagnetic quantities in $Q^f$. Since the electrical conductivity $\sigma$ satisfies $\sigma = \sigma^s = 0$ in $Q^s$, it also shows that the electromagnetic fields are not affected by the solid velocities $u^{s^i}$. In the boundary and interface conditions \eqref{831}, \eqref{832} the assumption of the magnetic permeability $\mu > 0$ being constant in the whole domain $Q$ (c.f.\ Section \ref{mechanicalsubsystem}) comes into play. Indeed, while the boundary condition for $B$ in \eqref{831} as well as the conditions \eqref{832} on $E$ are standard, the continuity of $B$ across the interface in \eqref{831} is not. However, it is standard to assume continuity of the normal component of $B$ and of the tangential component of $\frac{1}{\mu} B$ and hence, when $\mu$ is constant across $\partial S(t)$, we infer also the second relation in \eqref{831}. The reason why we cannot allow $B$ to have a jump across the interface is because otherwise it could not be a Sobolev function in $\Omega$, c.f.\ (\ref{-212}) below. \\[1em]
Both the mechanical subsystem \eqref{526}--\eqref{833} and the electromagnetic subsystem \eqref{94}--\eqref{832} can be studied in their own right, if one considers $B$ and $u$ as a prescribed external forcing term in Ohm's law \eqref{877} and the momentum equation \eqref{527}, respectively. In this article, however, we study the combined system \eqref{526}--\eqref{832}, in which both $B$ and $u$ are regarded as unknowns, coupled via Ohm's law and the presence of the Lorentz force $\frac{1}{\mu} \operatorname{curl}B \times B$ in the momentum equation. \\[1em]

{\centering \section{Variational formulation and main result} \label{weakformmainresult} \par }

\subsection{Notation} \label{notation}

The initial positions of the solid bodies are characterized through sets $S_0^i\subset \Omega $, $i=1,...,m$, onto which we impose the conditions
\begin{align}
\emptyset \neq S_0^i \ \text{is open and connected,}\quad |\partial S_0^i| = 0,\quad \overline{S}_0^i \bigcap \overline{S}_0^j = \emptyset \quad \quad \forall i,j=1,...,m,\ i \neq j. \label{-326}
\end{align}

Since the motion of the bodies is rigid, we can associate to each body an isometry $X^i(t,\cdot): \mathbb{R}^3 \rightarrow \mathbb{R}^3$, $t \in [0,T]$, such that its position $S^i(t)$ at an arbitrary time $t \in [0,T]$ can be expressed through the set-valued function
\begin{align}
S^i:\ [0,T] \rightarrow 2^{\mathbb{R}^3},\quad S^i(t) := X^i\left(t; S_0^i \right). \nonumber
\end{align}

In particular, with the notation
\begin{align}
X(t;\cdot):\ S_0 := \bigcup _{i=1}^m S_0^i \rightarrow \mathbb{R}^3,\quad X(t;\cdot)|_{S^i_0} := X^i(t;\cdot) \quad \quad \forall i = 1,...,m,\ t \in [0,T], \nonumber
\end{align}

the solid region at the time $t$ is given by $S(t):= X(t;S_0)$. In order to connect the motion of the bodies to the velocity field $u$ we require $u$ to be \textit{compatible} with the system $\lbrace S_0^i, X^i \rbrace_{i=1}^m$, i.e. we require the existence of rigid velocity fields $u^{s^i}(t,\cdot)$, $i=1,...,m$, such that
\begin{align}
u(t,x) = u^{s^i}(t,x) \quad \text{for a.a. } t \in [0,T]\ \text{and a.a. } x \in S^i(t) \label{-344}
\end{align}

and $X^i$ is the unique Carathéodory solution (c.f.\ \cite[Theorem 1.45]{roubicek}) to the initial value problem
\begin{align}
\frac{dX^i(t;x)}{dt} = u^{s^i}\left( t; X^i(t;x) \right), \quad X^i(0;x) = x,\quad x \in \mathbb{R}^3. \label{-345}
\end{align}

Finally we denote by $T(S)$ and $Y(S)$ the test function spaces
\begin{align}
T(S) :=& \left\lbrace \phi \in \mathcal{D}(Q):\ \mathbb{D}(\phi) = 0\ \text{in a neighbourhood of } \overline{Q}^s(S) \right\rbrace, \label{-339} \\
Y(S) :=& \left\lbrace b \in \mathcal{D}(Q):\ \operatorname{curl}b = 0\ \text{in a neighbourhood of } \overline{Q}^s(S) \right\rbrace \label{-340}
\end{align}

for our variational formulations of the momentum equation and the induction equation, respectively, in Definition \ref{weaksolutions} below.

\subsection{Weak solutions} \label{weaksolutionssection}

We are now in the position to present our variational formulation of the system \eqref{526}--\eqref{832}. With a slight abuse of notation we will write here and in the following sections $\sigma = \sigma^f > 0$, since the quantities containing $\sigma^s = 0$ are not visible in this weak formulation due to the non-conductivity of the solids.

\begin{definition}
\label{weaksolutions}

Let $T > 0$, let $\Omega \subset \mathbb{R}^3$ be a bounded domain and let $S_0 = \bigcup_{i=1}^m S_0^i$, where $S_0^i \subset \Omega$ for $i=1,...,m \in \mathbb{N}$ satisfy the conditions \eqref{-326}. Assume $\nu,\lambda, a, \gamma,\sigma,\mu \in \mathbb{R}$ to satisfy
\begin{align}
\nu, a, \sigma,\mu > 0,\quad \nu + \lambda \geq 0, \quad \gamma > \frac{3}{2}. \label{-324}
\end{align}

Moreover, consider some external forces $g,J \in L^\infty ((0,T)\times \Omega)$ and some initial data $0 \leq \rho_0 \in L^\gamma(\Omega)$, $(\rho u)_0 \in L^1(\Omega)$, $B_0 \in L^2(\Omega)$ such that
\begin{align}
\frac{\left|(\rho u)_0\right|^2}{\rho_0} \in L^1(\Omega),\quad (\rho u)_0 = 0 \ \ \text{a.e. in } \left\lbrace x \in \Omega:\ \rho_0(x) = 0 \right\rbrace,\quad \operatorname{div}B_0 = 0 \ \ \text{in } \mathcal{D}'(\Omega),\quad \left. B_0 \cdot \text{n}\right|_{\partial \Omega} = 0. \label{-325}
\end{align}

Then the system \eqref{526}--\eqref{832} is said to admit a weak solution if there exists a function
\begin{align}
X:\ [0,T] \times S_0 \rightarrow \mathbb{R}^3,\quad \left. X(t;\cdot) \right|_{S_0^i} = X^i(t;\cdot) \quad \quad \forall i=1,...,m,\ t\in [0,T], \label{-159}
\end{align}

where each $X^i(t;\cdot): \mathbb{R}^3 \rightarrow \mathbb{R}^3$ denotes an isometry, and if there exist functions
\begin{align}
0 \leq \rho &\in L^\infty \left( 0,T;L^\gamma (\Omega;\mathbb{R}) \right) \bigcap C\left([0,T];L^{1}(\Omega;\mathbb{R}) \right), \label{-211} \\
u &\in \left\lbrace \phi \in L^2\left(0,T;H_0^{1,2}\left(\Omega;\mathbb{R}^3 \right)\right):\ \mathbb{D}(\phi) = 0 \ \text{in } Q^s(S) \right\rbrace, \label{-213} \\
B &\in \left\lbrace b \in L^\infty \left(0,T;L^2\left(\Omega; \mathbb{R}^3\right)\right) \bigcap L^2\left(0,T;H_{\operatorname{div}}^{1,2}\left(\Omega;\mathbb{R}^3\right)\right):\ \operatorname{curl}b = 0 \ \text{in } Q^s(S),\ b \cdot \operatorname{n}|_{\partial \Omega} = 0 \right\rbrace, \label{-212}
\end{align}

where $S = S(\cdot) = X(\cdot, S_0)$, such that $\rho$ and $u$, extended by $0$ in $\mathbb{R}^3 \setminus \Omega$, satisfy the continuity equation and its renormalized form,
\begin{align}
\partial_t \rho + \operatorname{div} \left( \rho u \right) = 0 \quad &\text{in } \mathcal{D}'\left((0,T) \times \mathbb{R}^3 \right), \label{-214} \\
\partial_t \zeta (\rho) + \operatorname{div} \left( \zeta \left(\rho \right)u \right) + \left[ \zeta'\left(\rho \right)\rho - \zeta \left(\rho \right) \right] \operatorname{div} u = 0 \quad &\text{in } \mathcal{D}' \left((0,T) \times \mathbb{R}^3 \right), \label{-215}
\end{align}

for any
\begin{align}
\zeta \in C^1\left( [0,\infty) \right):\quad \left|\zeta'(r)\right| \leq cr^{\lambda _1} \quad \forall r \geq 1, \quad \quad \text{where } c>0,\ \lambda_1 > -1, \label{-216}
\end{align}

such that the momentum equation and the induction equation,
\begin{align}
- \int_0^T \int_\Omega \rho u \cdot \partial_t \phi \ dxdt =& \int _0^T \int _\Omega \left( \rho u \otimes u \right): \mathbb{D} (\phi) + a\rho^\gamma \operatorname{div} \phi - 2\nu \mathbb{D}(u):\mathbb{D} (\phi) \nonumber \\
&- \lambda \operatorname{div} u \operatorname{div} \phi + \rho g\cdot \phi + \frac{1}{\mu}\left( \operatorname{curl}B \times B \right) \cdot \phi \ dxdt, \label{-217} \\
- \int _0^T \int_\Omega B \cdot \partial_t b\ dxdt =& \int _0^T \int _\Omega \left[ -\frac{1}{\sigma \mu} \operatorname{curl} B + u \times B + \frac{1}{\sigma} J \right] \cdot \operatorname{curl} b \ dxdt, \label{-218}
\end{align}

are satisfied for any $\phi \in T(S)$ and any $b\in Y(S)$, such that the initial conditions
\begin{align}
\rho(0) = \rho_0,\quad \quad \left( \rho u \right)(0) = \left( \rho u \right)_0,\quad \quad B(0) = B_0, \label{-219}
\end{align}

hold true, where the latter two equations are to be understood in the sense that
\begin{align}
\lim_{t \rightarrow 0+} \int_\Omega \rho(t,x)u(t,x) \cdot \phi(x) \ dx = \int_\Omega (\rho u)_0(x) \cdot \phi (x)\ dx,\quad \lim_{t \rightarrow 0+} \int_\Omega B(t,x) \cdot b(x) \ dx = \int_\Omega B_0(x) \cdot b (x)\ dx \label{-341}
\end{align}

for all $\phi, b \in \mathcal{D}(\Omega)$ with $\mathbb{D}(\phi) = 0$ and $\operatorname{curl}b = 0$ in a neighbourhood of $\overline{S}_0$, respectively, and, finally, such that the system $\lbrace S_0^i, X^i \rbrace_{i=1}^m$ is compatible with the velocity field $u$.
\end{definition}

In this definition the compatibility of the velocity field $u$ and the system of rigid bodies leads to some vivid consequences for the solids. First of all, while the bodies are able to touch each other or the domain boundary, the possibility of interpenetrations are ruled out, c.f. \cite[Lemma 3.1, Corollary 3.1]{feireisl}. Moreover, even though the density does not satisfy a transport equation in the case of a compressible fluid, it still travels along the characteristics of $u$ in the solid part of the domain, c.f.\ \cite[Lemma 3.2]{feireisl}. For definiteness we present these results in the following lemma.

\begin{satz}[\cite{feireisl}] \label{compatibility}

Let $\Omega \subset \mathbb{R}^3$ and $S_0^i \subset \mathbb{R}^3$, $i=1,...,m \in \mathbb{N}$ be bounded domains and let further $u \in L^2(0,T;H_0^{1,2}(\Omega))$ be extended by $0$ outside of $\Omega$. Moreover, assume $u$ to be compatible with the system $\lbrace S_0^i, X^i \rbrace_{i=1}^m$ where each $X^i(t):\mathbb{R}^3 \rightarrow \mathbb{R}^3$, $t \in [0,T]$, $i=1,...,m$, denotes an isometry. Then it holds:
\begin{itemize}
\item[(i)] If, for $i \neq j \in \lbrace 1,...,m \rbrace$, there exists $\tau \in [0,T]$ such that $S^i(\tau) \bigcap S^j(\tau) \neq \emptyset$ then $X^i(t)= X^j(t)$ for all $t \in [0,T]$. Further, if there exists $\tau \in [0,T]$ such that $S^i(\tau) \not\subset \Omega$, then $X^i(t) = \operatorname{Id}$ for all $t \in [0,T]$.
\item[(ii)] If $\rho \in L^\infty(0,T;L^\gamma (\Omega))$, $\gamma > 1$, extended by $0$ outside of $\Omega$, satisfies
\begin{align}
\partial_t \rho + \operatorname{div} (\rho u) = 0 \quad \quad \text{in } \mathcal{D}'\left( (0,T) \times \mathbb{R}^3 \right) \nonumber
\end{align}
then
\begin{align}
\rho \left(t, X^i(t;x) \right) = \rho \left(0,x \right) \quad \quad \text{for all } t \in [0,T],\quad i=1,...,m \quad \text{and a.a. } x \in S_0^i. \label{-333}
\end{align}
\end{itemize}

\end{satz}

A detailed proof of the assertions (i) and (ii) is given in \cite[Lemma 3.1, Corollary 3.1]{feireisl} and \cite[Lemma 3.2]{feireisl}, respectively. The first part of assertion (i) can be shown directly from the fact that $\lbrace S_0^i, X^i \rbrace$ and $\lbrace S_0^j, X^j \rbrace$ are compatible with the same velocity field $u$, the second part then follows by regarding $\mathbb{R}^3\setminus \Omega$ as a rigid body with the associated rigid velocity field $0$. The proof of the assertion (ii) is achieved via a regularization of $\rho$ with respect to the spatial variable and a subsequent application of the regularization method by DiPerna and Lions, c.f.\ \cite{dipernalions}, to the continuity equation \eqref{-333} on compact subsets of the solid time-space domain.

\subsection{Main result} \label{mainresultsection}

Our main result yields existence of the weak solutions as introduced in Definition \ref{weaksolutions}.

\begin{theorem}
\label{mainresult}
Let $T > 0$, assume $\Omega \subset \mathbb{R}^3$ to be a simply connected domain of class $C^2 \bigcup C^{0,1}$ and assume $S_0^i \subset \Omega$, $i=1,...,m \in \mathbb{N}$ to be domains of class $C^2 \bigcup C^{0,1}$ which satisfy the conditions \eqref{-326}. Assume moreover the coefficients $\sigma, \mu, \nu, \lambda, a, \gamma \in \mathbb{R}$ to satisfy the conditions \eqref{-324} and the data $g,J \in L^\infty((0,T)\times \Omega)$, $\rho_0 \in L^\gamma (\Omega)$, $(\rho u)_0 \in L^1(\Omega)$ and $B_0 \in L^2(\Omega)$ to satisfy the conditions \eqref{-325}. Then the system \eqref{526}--\eqref{832} admits a weak solution $\lbrace X,\rho,u,B \rbrace$ in the sense of Definition \ref{weaksolutions} which in addition satisfies the energy inequality
\begin{align}
&\int_\Omega \frac{1}{2} \rho (t)\left|u (t)\right|^2 + \frac{a}{\gamma - 1} \rho ^\gamma(t) + \frac{1}{2\mu} \left| B(t) \right|^2\ dx + \int_0^t \int_\Omega 2\nu |\mathbb{D} \left(u \right)|^2 + \lambda \left| \operatorname{div}u \right|^2 + \frac{1}{\sigma \mu^2} \left| \operatorname{curl} B \right|^2\ dxd\tau \nonumber \\
\leq& \int_\Omega  \frac{1}{2} \frac{\left|\left(\rho u \right)_{0}\right|^2}{\rho_{0}} + \frac{a}{\gamma - 1} \rho_{0}^\gamma + \frac{1}{2\mu} \left| B_{0} \right|^2\ dx +\int_{0}^t\int_\Omega \rho g \cdot u + \frac{1}{\sigma \mu} J \cdot \operatorname{curl} B \ dxd\tau \label{-220}
\end{align}
 
for almost all $t \in [0,T]$.
\end{theorem}

The remainder of the article is devoted to the proof of Theorem \ref{mainresult}. In the following section we begin with an outline of the main ideas by introducing the approximation method on which the proof is based.

{\centering  \section{Approximate system} \label{approximatesystem} \par }

The biggest challenge in the extension of the proof in the incompressible case in \cite{incompressiblecase} to the compressible case lies in the correct construction of the approximate problem. We fix five parameters $\alpha, \epsilon, \eta > 0$, $n \in \mathbb{N}$ and $\Delta t > 0$ and introduce an approximation which consists of five different approximation levels, each of which corresponds to one of the parameters. The approximate system is chosen such that it is easy to solve; a solution to the original system is obtained by passing to the limit in all of the approximation levels. The first three approximation levels, associated to $\alpha, \epsilon$ and $\eta$, respectively, correspond to the approximation used in \cite{feireisl} for the purely mechanical problem: On the $\alpha$- and $\epsilon$-levels, the system is regularized through the addition of an artificial pressure term and multiple further regularization terms. The $\eta$-level consists of a penalization method which allows us to pass from a fluid-only system to a system containing both a fluid and rigid bodies. The fourth level, indexed by $n$, describes a Galerkin method used for solving the approximate momentum equation. Finally, on the fifth level, associated to $\Delta t$, the induction equation is discretized with respect to the time variable while the mechanical part of the problem is split up into a series of time-dependent problems on the small intervals between the discrete times. In the following we present the complete approximate system on the highest approximation level and subsequently give a more explicit description of each included approximation level and its purpose.

Let $\beta > 0$ be sufficiently large such that it satisfies in particular $\beta > \max\lbrace 4, \gamma \rbrace$ as in \cite[Section 6]{feireisl}. Let $V_n$, $n \in \mathbb{N}$, denote the Galerkin space spanned by the first $n$ eigenfunctions of the Lamé equation in $\Omega$ which constitute an orthonormal basis of $L^2(\Omega)$ and an orthogonal basis of $H_0^{1,2}(\Omega)$, c.f.\ \cite[Lemma 4.33]{novotnystraskraba}. Then, provided that the approximate system has already been solved up to the (discrete) time $(k-1)\Delta t$ for some $k = 1,...,\frac{T}{\Delta t}$, the approximate problem on the interval $[(k-1)\Delta t, k\Delta t]$ consists of finding a solution
\begin{align}
\rho_{\Delta t, k} &\in \left\lbrace \psi \in C\left([(k-1)\Delta t, k\Delta t]; C^{2,\frac{1}{2}} \left( \overline{\Omega} \right) \right) \bigcap C^1\left([(k-1)\Delta t, k\Delta t]; C^{0,\frac{1}{2}} \left( \overline{\Omega} \right) \right):\ \left. \nabla \psi \cdot \text{n} \right|_{\partial \Omega} = 0 \right\rbrace, \label{-317} \\
u_{\Delta t, k} &\in C\left([(k-1)\Delta t, k\Delta t]; V_n \right), \label{-316} \\
B_{\Delta t}^k &\in Y^k\left(S_{\Delta t} \right) :=\left\lbrace b \in H^{2,2}\left(\Omega \right):\ b \cdot \text{n}|_{\partial \Omega} = 0,\ \text{div}b=0,\ \text{curl} b = 0\ \text{in } S_{\Delta t}(k\Delta t) \bigcap \Omega \right\rbrace \label{-315}
\end{align}

to the system
\begin{align}
\partial_t \rho_{\Delta t, k} + \operatorname{div} \left( \rho_{\Delta t,k} u_{\Delta t,k} \right) =& \epsilon \Delta \rho_{\Delta t,k} \quad \text{in } [(k-1)\Delta t,k\Delta t] \times \Omega, \label{8} \\
\int_{(k-1)\Delta t}^{k\Delta t} \int_\Omega \partial_t \left(\rho_{\Delta t,k} u_{\Delta t,k} \right) \cdot \phi \ dxdt =& \int _{(k-1)\Delta t}^{k\Delta t} \int _\Omega \left( \rho_{\Delta t,k} u_{\Delta t,k} \otimes u_{\Delta t,k} \right): \mathbb{D} (\phi) + \left( a\rho_{\Delta t,k}^\gamma + \alpha \rho_{\Delta t,k}^\beta \right) \operatorname{div} \phi \nonumber \\
&- 2\nu \left( \chi_{\Delta t}^{k-1} \right)\mathbb{D}\left(u_{\Delta t,k}\right):\mathbb{D} (\phi) - \lambda \left( \chi_{\Delta t}^{k-1} \right)\operatorname{div}\left(u_{\Delta t,k}\right)\operatorname{div} \phi \nonumber \\
&+ \rho_{\Delta t,k} g\cdot \phi + \frac{1}{\mu}\left( \text{curl}B_{\Delta t}^{k-1} \times B_{\Delta t}^{k-1} \right) \cdot \phi - \epsilon \left| u_{\Delta t,k} \right|^2u_{\Delta t,k} \cdot \phi \nonumber \\
&- \epsilon \left( \nabla u_{\Delta t,k} \nabla \rho_{\Delta t,k} \right) \cdot \phi\ dxdt, \label{57} \\
-\int_\Omega \frac{B_{\Delta t}^k - B_{\Delta t}^{k-1}}{\Delta t} \cdot b\ dx =& \int_\Omega \left[ \frac{1}{\sigma \mu} \operatorname{curl} B_{\Delta t}^k  - \tilde{u}_{\Delta t}^{k-1} \times B_{\Delta t}^{k-1} \right. \nonumber \\
&+ \left. \frac{\epsilon}{\mu ^2} \left| \operatorname{curl} B_{\Delta t}^k \right|^2\operatorname{curl} B_{\Delta t}^k -\frac{1}{\sigma} J_{\Delta t}^k \right] \cdot \operatorname{curl} b \nonumber \\
&+ \epsilon \operatorname{curl}\left(\operatorname{curl} B_{\Delta t}^k\right) \cdot \operatorname{curl}\left(\operatorname{curl} b\right)\ dx \label{116}
\end{align}

for all $\phi \in C([(k-1)\Delta t, k\Delta t];V_n)$ and
\begin{align}
b \in W^k\left(S_{\Delta t}\right) :=& \bigg \lbrace b\in H^{2,2}\left(\Omega \right):\ b \cdot \text{n}|_{\partial \Omega} = 0,\ \text{curl}\ b = 0\ \text{in } S_{\Delta t}(k\Delta t) \bigcap \Omega \bigg \rbrace, \label{937}
\end{align}

which in addition satisfies the initial conditions
\begin{align}
\rho_{\Delta t,k}((k-1)\Delta t;x) &= \rho_{\Delta t,k-1}((k-1)\Delta t; x),\quad \rho_{\Delta t,1}(0;x) = \rho_0(x),\quad x \in \Omega, \label{-319} \\
u_{\Delta t,k}((k-1)\Delta t;x) &= u_{\Delta t,k-1}((k-1)\Delta t; x),\quad u_{\Delta t,1}(0;x) = u_0(x),\quad x \in \Omega, \label{-320} \\
B_{\Delta t}^0(x) &= B_0(x),\quad x \in \Omega. \label{-321}
\end{align}

Before we proceed with the explanation of the different approximation levels in \eqref{-317}--\eqref{-321}, let us clarify the notation introduced in this system: For the definition of the set $S_{\Delta t}(k\Delta t)$ in \eqref{-315} and \eqref{937}, we first denote by
\begin{align}
O^i := (S^i_0)_\delta := \left\lbrace x \in S^i_0:\ \operatorname{dist}\left(x, \partial S_0^i \right) > \delta \right\rbrace \nonumber
\end{align}

the $\delta$-kernel of the initial domain $S^i_0$ of the $i$-th body, where $\delta > 0$ is chosen sufficiently small such that for all $i=1,...,m$ the $\delta$-neighbourhood $(O^i)^\delta = \lbrace x \in \mathbb{R}^3:\ \operatorname{dist}(x, (O^i)^\delta) < \delta \rbrace$ of $O^i$ coincides with $S^i_0$. Such $\delta > 0$ exists due to the $C^2$-regularity of $S^i_0$, c.f. \cite[Proposition 2.1]{tucsnak}. Then we set
\begin{align}
O := \bigcup_{i=1}^m O^i. \nonumber
\end{align}

Moreover, we denote by $X_{\Delta t,k}$ the unique solution to the initial value problem
\begin{align}
\frac{d}{dt}X_{\Delta t,k}(t;x) =& R_\delta \left[ u_{\Delta t, k} \right] \left( t, X_{\Delta t,k}(t;x) \right),\quad t \in [(k-1)\Delta t,k\Delta t], \label{1} \\
X_{\Delta t,k}((k-1)\Delta t;x) &= X_{\Delta t,k-1}((k-1)\Delta t; x),\quad X_{\Delta t,1}(0;x) = x,\quad x \in \mathbb{R}^3, \label{243}
\end{align}

where $R_\delta [u_{\Delta t,k}](t,\cdot) := u_{\Delta t,k}(t,\cdot)*\Theta_\delta (\cdot)$ and $\Theta_\delta$ denotes a radially symmetric and non-increasing mollifier with respect to the spatial variable. With this notation at hand we define the domain $S^i_{\Delta t}(t)$ of the $i$-th approximate solid at an arbitrary time $t \in [(k-1)\Delta t, k \Delta t] \subset [0,T]$ by
\begin{align}
O_{\Delta t}^i (t) := X_{\Delta t,k}\left(t;O^i\right),\quad S^i_{\Delta t}(t) := \left( O_{\Delta t}^i (t) \right)^\delta = \left\lbrace x \in \mathbb{R}^3:\ \operatorname{dist}\left(x, O^i_{\Delta t}(t) \right) < \delta \right\rbrace. \label{-323}
\end{align}

Consequently, the approximate solid region at time $t \in [0,T]$ is given by
\begin{align}
S_{\Delta t}(t) := \bigcup_{i=1}^m S^i_{\Delta t}(t), \nonumber
\end{align}

which in particular defines the set $S_{\Delta t}(k\Delta t)$ in \eqref{-315} and \eqref{937}. We note that, by construction, $S_{\Delta t}(t)$ can be an arbitrary subset of $\mathbb{R}^3$, while the corresponding approximate solid time-space domain $Q^s(S_{\Delta t})$, defined according to \eqref{-353}, is always a subset of the bounded domain $Q$. For later use we further remark that $S_{\Delta t}^i(t)$, as the $\delta$-neighbourhood of a bounded set, satisfies the cone condition and thus has the property
\begin{align}
\left| \partial S^i_{\Delta t}(t) \right| = 0 \quad \text{for all } t \in [0,T],\ i=1,...,m. \label{-318}
\end{align}

Next, for the definition of the variable viscosity coefficients $\nu ( \chi_{\Delta t}^{k-1})$ and $\lambda ( \chi_{\Delta t}^{k-1})$ in the momentum equation \eqref{57} we denote the signed distance function of arbitrary sets $U \subset \mathbb{R}^3$ by
\begin{align}
\textbf{db}_U(x) := \operatorname{dist}\left( x, \overline{\mathbb{R}^3 \setminus U} \right) - \operatorname{dist}\left( x, \overline{U} \right). \nonumber
\end{align}

Further we introduce the signed distance function of the approximate solid area,
\begin{align}
\chi_{\Delta t}(t,x) := \textbf{db}_{S_{\Delta t}(t)}(x),\quad \quad \chi_{\Delta t}^{k-1}(t) := \chi_{\Delta t}((k-1)\Delta t,\cdot) \quad \text{for } t \in [(k-1)\Delta t,k\Delta t]. \nonumber
\end{align}

Choosing a convex function $H \in C^\infty (\mathbb{R})$ such that
\begin{align}
H(z) = 0 \quad \text{for } z \in (-\infty, 0],\quad \quad H(z) > 0 \quad \text{for } z \in (0,+\infty), \label{-322}
\end{align}

we then define the variable viscosity coefficients by
\begin{align}
\nu \left( \chi_{\Delta t}^{k-1} \right) := \nu + \frac{1}{\eta} H\left(\chi_{\Delta t}^{k-1} \right),\quad \quad \lambda \left( \chi_{\Delta t}^{k-1} \right) := \lambda + \frac{1}{\eta} H\left(\chi_{\Delta t}^{k-1}\right). \label{-52}
\end{align}

Finally, in the induction equation \eqref{116} the function $\tilde{u}_{\Delta t}^{k-1}$ is defined by
\begin{align}
\tilde{u}_{\Delta t}^{k-1}(x) := \left\{
                \begin{array}{ll}
                  \frac{1}{\Delta t} \displaystyle\int_{(k-2)\Delta t}^{(k-1)\Delta t} u_{\Delta t, k-1}(\tau,x) d\tau \ \ \ &\text{if } k \geq 2,\\[5mm]
                  u_{\Delta t}^0(x) &\text{if } k = 1.
                \end{array}
              \right. \label{-235}
\end{align}

while the discretized external force $J_{\Delta t}^k$ is defined by
\begin{align}
J_{\Delta t}^k := J_{\omega}(k\Delta t), \quad J_\omega (t):= \int_0^T \theta _\omega \left( t + \omega \frac{T-2t}{T} - s \right)J(s)\ ds, \label{840}
\end{align}

for another mollifier $\theta_\omega :\mathbb{R}\to \mathbb{R}$ and a suitable choice of $\omega = \omega ( \Delta t),$ $\omega (\Delta t)\rightarrow 0$ for $\Delta t \rightarrow 0$. We are now in the position to discuss the several approximation levels and the reasons why they are required. We start from the highest level. \\
The $\Delta t$-level constitutes the level which contains most of the difficulties. It is here where the main novelties of our proof enter, compared to the incompressible setting in \cite{incompressiblecase}. The situation presents itself in the following way: On the one hand, the dependence of the test functions \eqref{-339}, \eqref{-340} for the induction equation and the momentum equation on the solution of the system hinders the effort to solve all of the equations in the system simultaneously. While we can deal with the test functions in the momentum equation by means of a penalization method (c.f.\ the $\eta$-level below), the same does not work in case of the induction equation. This suggests to decouple the system by the use of a classical time discretization, e.g.\ via the Rothe method (\cite[Section 8.2]{roubicek}). In this way, at each fixed discrete time we can first determine a velocity field and, from this, the position of the approximate solid. This in turn determines the test functions \eqref{937} and solving the discretized induction equation \eqref{116} becomes a routine matter. On the other hand, however, the various functions evaluated at different discrete times in a fully discretized system complicate the derivation of a meaningful energy inequality. The author could not find a way to transfer several of the techniques known for the continuous compressible Navier-Stokes system (c.f.\ \cite[Sections 7.6.5, 7.6.6, 7.7.4.2]{novotnystraskraba}) - in particular, the proof of the non-negativity of the density - to the discrete case and it did not seem to be possible to derive the uniform bounds required for the limit passage with respect to $\Delta t \rightarrow 0$. \\
Our solution to this dilemma consists of considering, instead of a strictly discretized system, a hybrid system in which the induction equation \eqref{116} is indeed discretized by the Rothe method, while the continuity equation and the momentum equation are solved as continuous equations on the small intervals between each pair of consecutive discrete times, c.f \eqref{8} and \eqref{57}. Through this, the solution dependence of the test functions in the induction equation can be handled as in the fully discrete system, while the mechanical part of the energy inequality - with the density bounded away from zero - can be derived as in the strictly continuous case. Moreover, under the consideration of piecewise linear interpolants of the discrete functions, the discrete induction equation \eqref{116} also leads to a continuous energy estimate, which can be combined with the mechanical estimate to obtain the full energy inequality, c.f.\ Section \ref{discreteenergyinequality}. A hybrid approximation scheme was already used in our proof in the incompressible case \cite{incompressiblecase}. In that case, however, the major part of the system could be discretized in time while only the transport equation for the characteristic function of the solid domain had to be treated as a continuous problem on small time intervals. The idea for the latter procedure, in turn, stems from \cite{gigli}. \\
The Galerkin method carried out on the $n$-level is used to solve the continuous momentum equation \eqref{8} on the small time intervals from the $\Delta t$-level by a standard procedure. The Galerkin-regularity of the velocity field furthermore helps us during the limit passage with respect to $\Delta t \rightarrow 0$, c.f.\ \eqref{-259} below. \\
After letting $n$ tend to $\infty$ we find ourselves in the same situation as in the approximation of the exclusively mechanical system in \cite[Section 6]{feireisl}. Indeed, the remaining three approximation levels correspond directly to the three level approximation scheme used in that article. Hence, for the mechanical part of our problem we can follow exactly the strategy used therein. Moreover, the limit passages in the induction equation from here on do not contain any new difficulties anymore. Consequently, after the limit passage in $n$ the rest of the proof will become a routine matter. \\
The penalization method on the $\eta$-level - c.f.\ Section \ref{etalimit} - is the same as the one used for the fluid-rigid bodies system in \cite{feireisl} and was, before that, also used for example for the corresponding two dimensional problem in \cite{tucsnak}. The idea behind it is to approximate the entirety of the fluid and the rigid bodies by a fluid in the whole domain with viscosity tending to infinity in the later solid regions. Mathematically this is implemented through the variable viscosity coefficients \eqref{-52}. Due to the choice of the function $H$ in \eqref{-322} these coefficients blow up in the approximate solid region once we let $\eta$ tend to $0$ and, thanks to the energy inequality, this will cause the limit velocity field $u$ to coincide with a rigid velocity field in each body. Moreover, the positions $S^i(t)$ of the bodies in the $\eta$-limit are determined through the flow curves of $R_\delta [u]$, c.f.\ \eqref{1} and \eqref{-323}. This regularized velocity field has the useful property that, for any domain $U \subset \mathbb{R}^3$, it holds
\begin{align}
\mathbb{D}(u(t,\cdot)) = 0\quad \text{in } U \quad \quad \Rightarrow \quad \quad R_\delta [u](t,\cdot) = u(t,\cdot) \quad \text{in } U_\delta = \left\lbrace x \in U:\ \operatorname{dist}\left( x, \partial U \right)>\delta \right\rbrace, \label{-348}
\end{align}

c.f.\ \cite[Remark 6.1]{feireisl}. Hence $R_\delta[u]$ coincides with $u$ itself in the sets $O^i(t) \subset S^i(t)$, in which $\mathbb{D}(u)=0$. Consequently, the rigid velocity fields coinciding with $u$ in $S^i(t)$ also coincide, in $O^i(t)$, with the velocity field $R_\delta [u]$ which determines the motion of the bodies. In particular, this shows that the bodies $S^i(t)$ are indeed rigid. \\
On the $\epsilon$-level, the continuity equation \eqref{8} is regularized through the additional Laplacian $\epsilon \Delta \rho_{\Delta t,k}$. The additional quantity $\epsilon ( \nabla u_{\Delta t,k} \nabla \rho_{\Delta t,k})$ in the momentum equation \eqref{57} ensures that the energy inequality is preserved under this regularization. This procedure, which is classical in the theory of the compressible Navier-Stokes equations, is what guarantees us the non-negativity of the density, c.f.\ \cite[Section 7.3.8]{novotnystraskraba}. The other regularization term $\epsilon |u_{\Delta t,k}|^2u_{\Delta t,k}$ in \eqref{57} is needed in the time discrete level where, as opposed to the continuous case, the mixed terms from the momentum equation and the induction equation do not annihilate each other in the energy inequality, which prevents a direct application of the Gronwall lemma. The quantity $\epsilon |u_{\Delta t,k}|^2u_{\Delta t,k}$ can be used to control the velocity part of these mixed terms. The $4$-double-curl $\epsilon\operatorname{curl}(\left| \operatorname{curl} B_{\Delta t}^k \right|^2\operatorname{curl} B_{\Delta t}^k)$ in the induction equation \eqref{116} fulfills, as in the incompressible setting in \cite{incompressiblecase}, the same purpose for the magnetic part of the mixed terms so that we are able to derive uniform bounds from the energy inequality nevertheless, c.f.\ Section~\ref{discreteenergyinequality}. We remark that this control of the mixed terms was also the motivation for the definition of the velocity field \eqref{-235} in the discrete induction equation: Indeed, defining this quantity as a mean value of the velocity field obtained from the momentum equation on the intervals $[(k-2)\Delta t, (k-1)\Delta t]$, we can absorb it into the left-hand side of the energy inequality thanks to the above-mentioned regularization terms. If instead the term was defined, more intuitively, as a pointwise evaluation of $u_{\Delta t,k}$, we would not be able to handle it. The last regularization term in \eqref{116}, the $4$-th curl of $B_{\Delta t}^k$, enables us to express the induction equation via some weakly continuous operator on $Y^k(S_{\Delta t})$. Seeing that this operator is moreover coercive, we will then be able to infer the existence of $B_{\Delta t}^k$, c.f.\ Section~\ref{existenceB}. \\
Finally, on the $\alpha$-level, the artificial pressure term $\alpha \rho_{\Delta t,k}^\beta$ is added to the momentum equation \eqref{57}. Again this method is already well-known from the general existence theory for the compressible Navier-Stokes system, c.f.\ \cite[Section 7.3.8]{novotnystraskraba}. The artificial pressure gives us an additional amount of integrability of the density and its gradient, required to pass to the limit in the term $\epsilon ( \nabla u_{\Delta t,k} \nabla \rho_{\Delta t,k})$ from the $\epsilon$-level, c.f.\ \cite[Section 7.8.2]{novotnystraskraba}. It furthermore simplifies the limit passage with respect to $\epsilon \rightarrow 0$, since the additional integrability allows for the use of the regularization technique by DiPerna and Lions, c.f.\ \cite[Lemma 6.8, Lemma 6.9]{novotnystraskraba}. \\

{\centering \section{Existence of the approximate solution} \label{existencediscretesolution} \par }

We begin the proof of Theorem \ref{mainresult} by showing the existence of a solution to the approximate problem \eqref{-317}--\eqref{-321} on the highest approximation level. \\

\subsection{Existence of the density and velocity}

The existence of the density and the velocity field on the Galerkin level can be shown by classical methods, c.f.\ for example \cite[Section 7.7]{novotnystraskraba}. More precisely, the continuity equation \eqref{8} and the momentum equation \eqref{57} can be solved simultaneously by means of a fixed point argument: For fixed $w \in C([(k-1)\Delta t,k\Delta t];V_n)$ we consider the Neumann problem
\begin{align}
&\ \quad \quad \quad \partial_t \rho + \operatorname{div} \left( \rho w \right) = \epsilon \Delta \rho \quad \quad \text{in } [(k-1)\Delta t,k\Delta t] \times \Omega, \label{-228} \\
&\left. \nabla \rho \cdot \text{n}\right|_{\partial \Omega} = 0,\quad \rho((k-1)\Delta t, \cdot) = \rho_{\Delta t,k-1}((k-1)\Delta t; \cdot)\quad \quad \text{in } \Omega, \label{-352} \\
&\ \quad \quad \quad 0 < \underline{\rho} \leq \rho_{\Delta t,k-1}((k-1)\Delta t; x) \leq \overline{\rho} < \infty, \quad \quad \text{in } \Omega. \label{-231}
\end{align}

It is well known (c.f.\ Lemma \ref{neumannproblem} in the appendix) that \eqref{-228}--\eqref{-231} admits a unique solution
\begin{align}
\rho = \rho (w) \in C\left([(k-1)\Delta t, k\Delta t]; C^{2,\frac{1}{2}} \left( \overline{\Omega} \right) \right) \bigcap C^1\left([(k-1)\Delta t, k\Delta t]; C^{0,\frac{1}{2}} \left( \overline{\Omega} \right) \right), \nonumber
\end{align}

which satisfies the estimate
\begin{align}
0 < \underline{\rho} \exp \left(-\int_{(k-1)\Delta t}^t \left\| w(\tau) \right\|_{V_n}\ d\tau \right) \leq \rho(w)(t,x) \leq \overline{\rho} \exp \left(\int_{(k-1)\Delta t}^t \left\| w(\tau) \right\|_{V_n}\ d\tau \right) < \infty \label{-233}
\end{align}

for all $(t,x) \in [0,T]\times \overline{\Omega}$. Further, we consider a linearized version of the momentum equation \eqref{57}. Given $w \in C([(k-1)\Delta t,k\Delta t];V_n)$ and the associated solution $\rho (w)$ to the Neumann problem \eqref{-228}--\eqref{-231}, we seek $u \in C([(k-1)\Delta t,k\Delta t];V_n)$ such that
\begin{align}
\int_\Omega \partial_t \left(\rho(w) u \right) \cdot \phi \ dx =& \int _\Omega \left( \rho(w) w \otimes u \right): \mathbb{D} (\phi) + \left(a\rho^\gamma (w) + \alpha \rho^\beta (w)\right) \operatorname{div} \phi \nonumber \\
& - 2\nu \left( \chi_{\Delta t}^{k-1} \right)\mathbb{D}\left(u\right):\mathbb{D} (\phi) - \lambda \left( \chi_{\Delta t}^{k-1} \right)\operatorname{div}\left(u\right)\operatorname{div} \phi \nonumber \\
&+ \rho (w) g\cdot \phi + \frac{1}{\mu}\left( \text{curl}B_{\Delta t}^{k-1} \times B_{\Delta t}^{k-1} \right) \cdot \phi \nonumber \\
&- \epsilon \left( \nabla u \nabla \rho(w) \right) \cdot \phi - \epsilon \left| w \right|^2u \cdot \phi\ dx \quad \quad \text{in } [(k-1)\Delta t,k\Delta t], \nonumber \\
u((k-1)\Delta t,\cdot) =& u_{\Delta t,k-1}((k-1)\Delta t; \cdot) \quad \quad \quad \quad \quad \quad \quad \quad \ \text{in } \Omega, \nonumber
\end{align}

for all $\phi \in C([(k-1)\Delta t,k\Delta t];V_n)$. Under exploitation of the fact that, by \eqref{-233}, $\rho(w)$ is bounded away from $0$ and the linearity of the problem, it follows from classical methods that this problem admits a unique solution $u = u(w) \in C([(k-1)\Delta t,k\Delta t];V_n)$. We can thus define an operator
\begin{align}
\mathbb{T}: C\left([(k-1)\Delta t, k \Delta t];V_n\right) \rightarrow C\left([(k-1)\Delta t, k \Delta t];V_n\right),\quad \mathbb{T}(w) := u. \nonumber
\end{align}

It is easy to see that $\mathbb{T}$ is continuous and compact and, by an energy estimate, fixed points of $s\mathbb{T}$ for $s \in [0,1]$ are bounded in $C([(k-1)\Delta t,k\Delta t];V_n)$, uniformly with respect to $s$. Under these conditions the Schaefer fixed point theorem (see \cite[Section 9.2.2, Theorem 4]{evans}) tells us that $\mathbb{T}$ possesses a fixed point $u_{\Delta t,k} \in C([(k-1)\Delta t,k\Delta t];V_n)$, which constitutes the desired solution to the initial value problem \eqref{57}, \eqref{-320}. Furthermore, by construction, the associated density $\rho_{\Delta t,k} := \rho (u_{\Delta t,k})$ is the desired solution to the corresponding initial value problem \eqref{8}, \eqref{-319} for the density.

\subsection{Existence of the magnetic induction} \label{existenceB}

The existence of the magnetic induction is obtained as in the incompressible case, c.f.\ \cite[Section 3]{incompressiblecase}. We equip the space $Y^k(S_{\Delta t})$ with the norm $\|\cdot \|_{H^{2,2}(\Omega)}$ and express the identity \eqref{116} through the equation
\begin{equation}
\left\langle AB_{\Delta t}^k, b \right\rangle_{\left(Y^k(S_{\Delta t})\right)^* \times Y^k(S_{\Delta t})} = \left\langle f, b \right\rangle_{\left(Y^k(S_{\Delta t})\right)^* \times Y^k(S_{\Delta t})} \quad \forall b \in Y^k\left(S_{\Delta t}\right), \label{109}
\end{equation}

where the operator $A:Y^k(S_{\Delta t}) \rightarrow (Y^k(S_{\Delta t}))^*$ and the right-hand side $f \in (Y^k(S_{\Delta t}))^*$ are defined by
\begin{align}
\left\langle A(B),b \right\rangle _{\left(Y^k(S_{\Delta t})\right)^* \times Y^k(S_{\Delta t})} :=& \int _\Omega  \frac{B}{\Delta t} \cdot b + \epsilon \text{curl}\left(\text{curl} B \right) \cdot \text{curl}\left(\text{curl} b\right) \nonumber \\
&+ \left[ \frac{1}{\sigma \mu} \text{curl}B + \frac{\epsilon}{\mu ^2} \left| \text{curl} B \right|^2\text{curl} B \right] \cdot \text{curl}b\ dx, \nonumber \\
\left\langle f,b \right\rangle _{\left(Y^k(S_{\Delta t})\right)^* \times Y^k(S_{\Delta t}) } :=& \int _\Omega \frac{B_{\Delta t}^{k-1}}{\Delta t} \cdot b +\left( \tilde{u}_{\Delta t}^{k-1} \times B_{\Delta t}^{k-1} \right) \cdot \text{curl}b + \frac{1}{\sigma}J_{\Delta t}^k \cdot \operatorname{curl}b\ dx \nonumber
\end{align}

for any $B,b \in Y^k(S_{\Delta t})$. The operator $A$ is weakly continuous and coercive and consequently surjective from $Y^k(S_{\Delta t})$ onto $(Y^k(S_{\Delta t}))^*$, c.f.\ \cite[Theorem 1.2]{francu}. In particular, there exists a function $B_{\Delta t}^k \in Y^k(S_{\Delta t})$ which satisfies \eqref{109} and thus the induction equation \eqref{116} for all $b \in Y^k(S_{\Delta t})$. Finally, by the Helmholtz decomposition, see \cite[Theorem 4.2]{maxwellinequality}, we infer that \eqref{116} does not only hold for $b \in Y^k(S_{\Delta t})$ but also for the (not divergence-free) test functions $b \in W^k(S_{\Delta t})$. Altogether, we have shown the following result:
\begin{proposition}
\label{System on Delta t-level}
Let $n \in \mathbb{N}$, $\Delta t, \eta, \epsilon, \alpha > 0$ such that $\frac{T}{\Delta t} \in \mathbb{N}$ and let $\beta > \max\lbrace 4, \gamma \rbrace$ be sufficiently large. Assume the conditions of Theorem \ref{mainresult} to be satisfied. Moreover, assume that
\begin{align}
\rho _0 \in& C^{2,\frac{1}{2}}\left(\overline{\Omega}\right),\quad \quad (\rho u)_0 \in C^2\left(\overline{\Omega}\right), \quad \quad u_0 := P_n \left(\frac{(\rho u)_0}{\rho_0}\right) \in V_n, \quad \quad B_0 \in H^{2,2}(\Omega), \nonumber \\
&0 < \alpha \leq \rho_0 \leq \alpha^{-\frac{1}{2\beta}},\quad \quad \left. \nabla \rho_0 \cdot \operatorname{n} \right|_{\partial \Omega} = 0, \quad \quad \operatorname{div} B_0 = 0,\quad \quad \left. B_0 \cdot \operatorname{n}\right|_{\partial \Omega} = 0, \nonumber
\end{align}

where $P_n$ denotes the orthogonal projection of $L^2(\Omega)$ onto $V_n$. Finally, for all $k = 1,...,\frac{T}{\Delta t}$, let $J_{\Delta t}^k$ be defined by \eqref{840}. Then, for each $k = 1,...,\frac{T}{\Delta t}$, there exist functions
\begin{align}
0 \leq \rho_{\Delta t, k} &\in \bigg\lbrace \psi \in C\left([(k-1)\Delta t, k\Delta t]; C^{2,\frac{1}{2}} \left( \overline{\Omega} \right) \right) \bigcap C^1\left([(k-1)\Delta t, k\Delta t]; C^{0,\frac{1}{2}} \left( \overline{\Omega} \right) \right):\ \left. \nabla \psi \cdot \operatorname{n} \right|_{\partial \Omega} = 0 \bigg\rbrace, \nonumber \\
&\quad \quad \quad u_{\Delta t, k} \in C\left([(k-1)\Delta t, k\Delta t]; V_n \right), \quad \quad B_{\Delta t}^k \in Y^k\left(S_{\Delta t} \right) \nonumber,
\end{align}

which satisfy the continuity equation \eqref{8}, the momentum equation \eqref{57} for all test functions $\phi \in C\left([(k-1)\Delta t, k\Delta t]; V_n \right)$ and the induction equation \eqref{116} for all test functions $b \in W^k(S_{\Delta t})$ as well as the initial conditions \eqref{-319}--\eqref{-321}.

\end{proposition}

{\centering \section{Limit passage in the time discretization} \label{deltatlimit} \par }

We continue by passing to the limit with respect to $\Delta t \rightarrow 0$. As in \cite{incompressiblecase}, we first need to assemble the functions constructed in Section \ref{existencediscretesolution}, defined up to now only on small time intervals or in discrete time points, to functions defined on the whole time interval $[0,T]$. More precisely, for functions $f_{\Delta t,k}$, defined on $[(k-1)\Delta t, k\Delta t] \times \Omega$ for $k=1,...,\frac{T}{\Delta t}$, we introduce the assembled functions
\begin{align}
f_{\Delta t} (t,\cdot) := f_{\Delta t,k}(t,\cdot) \quad \quad \forall t \in [(k-1)\Delta t, k\Delta t),\quad k=1,...,\frac{T}{\Delta t} \label{-245}
\end{align}

while for discrete functions $h_{\Delta t}^k$, defined on $\Omega$ for $k=0,...,\frac{T}{\Delta t}$, we introduce the piecewise affine and piecewise constant interpolants
\begin{align}
h_{\Delta t}(t) &:= \left( \frac{t}{\Delta t} - (k-1)\right)h_{\Delta t}^k + \left(k - \frac{t}{\Delta t}\right) h_{\Delta t}^{k-1}\ \ \ &\forall& t \in [(k-1)\Delta t, k\Delta t),\quad k=1,...,\frac{T}{\Delta t}, \label{21} \\
\overline{h}_{\Delta t}(t) &:= h_{\Delta t}^k \ \ \ &\forall& t \in [(k-1)\Delta t, k\Delta t),\quad k=0,...,\frac{T}{\Delta t}, \label{22} \\
\overline{h}'_{\Delta t}(t) &:= h_{\Delta t}^{k-1} \ \ \ &\forall& t \in [(k-1)\Delta t, k\Delta t),\quad k=1,...,\frac{T}{\Delta t}. \label{23}
\end{align}

Moreover, in order to derive a suitable energy inequality in Section \ref{discreteenergyinequality} below, we also introduce a piecewise affine interpolation of the square of the $L^2(\Omega)$-norm,
\begin{align}
h_{\Delta t, \|\cdot \|}(t) := \left( \frac{t}{\Delta t} - (k-1) \right) \left\| h_{\Delta t}^k \right\|_{L^2(\Omega)}^2 + \left( k - \frac{t}{\Delta t} \right) \left\| h_{\Delta t}^{k-1} \right\|_{L^2(\Omega)}^2,\quad \quad \forall t \in [(k-1)\Delta t, k\Delta t) \nonumber
\end{align}

for any $k=1,...,\frac{T}{\Delta t}.$ Since, by Proposition \ref{System on Delta t-level}, the functions $\rho_{\Delta t,k}$ and $u_{\Delta t,k}$ satisfy the continuity equation \eqref{8}, the momentum equation \eqref{57}, and the initial conditions \eqref{-319}--\eqref{-321}, it follows from the definition of $\rho_{\Delta t}$ and $u_{\Delta t}$ in \eqref{-245} as well as of $B_{\Delta t}$, $\overline{B}_{\Delta t}$ and $\overline{B}'_{\Delta t}$ in \eqref{21}--\eqref{23} that these functions solve the continuity equation
\begin{align}
\partial_t \rho_{\Delta t} + \operatorname{div} \left( \rho_{\Delta t} u_{\Delta t} \right) = \epsilon \rho_{\Delta t} \quad \text{a.e. in } [0,T] \times \Omega, \label{-246}
\end{align}

the momentum equation
\begin{align}
\int_0^T \int_\Omega \partial_t \left( \rho_{\Delta t} u_{\Delta t} \right) \cdot \phi \ dxdt =& \int _0^T \int _\Omega \left( \rho_{\Delta t} u_{\Delta t} \otimes u_{\Delta t} \right): \mathbb{D} (\phi) + \left( a\rho_{\Delta t}^\gamma + \alpha \rho_{\Delta t}^\beta \right) \operatorname{div} \phi \nonumber \\
&- 2\nu \left(\overline{\chi}_{\Delta t}' \right) \mathbb{D}(u_{\Delta t}):\mathbb{D} (\phi) - \lambda \left(\overline{\chi}_{\Delta t}' \right) \operatorname{div} u_{\Delta t} \operatorname{div} \phi - \epsilon \left| u_{\Delta t} \right|^2 u_{\Delta t} \cdot \phi \nonumber \\
&+ \rho_{\Delta t} g\cdot \phi + \frac{1}{\mu}\left( \operatorname{curl}\overline{B}_{\Delta t}' \times \overline{B}_{\Delta t}' \right) \cdot \phi - \epsilon \left( \nabla u_{\Delta t} \nabla \rho_{\Delta t} \right) \cdot \phi \ dxdt, \label{-248}
\end{align}

for any $\phi \in C([0,T];V_n)$ and the initial conditions
\begin{align}
\rho_{\Delta t}(0) = \rho_0,\quad u_{\Delta t}(0) = u_0,\quad B_{\Delta t}(0)=B_0. \label{-247}
\end{align}

Furthermore, from $X_{\Delta t,k}$ being the unique solution to the initial value problem \eqref{1}, \eqref{243}, it follows that $X_{\Delta t}$ is the unique solution to
\begin{align}
\frac{dX_{\Delta t}(t;x)}{dt} = R_\delta \left[ u_{\Delta t} \right] \left( t, X_{\Delta t}(t;x) \right),\quad \text{for } t \in [0,T], \quad \quad X_{\Delta t}(0;x) =  x,\quad \text{for } x \in \mathbb{R}^3. \label{-254}
\end{align}

Finally, we consider functions
\begin{align}
b \in L^4\left(0,T;H_0^{2,2}(\Omega) \right) \quad \text{such that} \quad b(\tau) \in W^l\left(S_{\Delta t} \right) \quad \text{for a.a. } \tau \in [(l-1)\Delta t,l\Delta t],\ l=1,...,\frac{T}{\Delta t} \label{-41}
\end{align}

and realize that, after a density argument, the discrete induction equation \eqref{116} at time $k\Delta t$ can be tested by $b(t)$ for almost all $t \in [(k-1)\Delta t, k\Delta t]$. After integration over $[(k-1)\Delta t, k\Delta t]$ and summation over $k$ this yields
\begin{align}
\int_0^T \int_\Omega \partial_t B_{\Delta t} \cdot b\ dxdt =& \int_0^T\int_\Omega \left[ -\frac{1}{\sigma \mu} \operatorname{curl} \overline{B}_{\Delta t} + \overline{\tilde{u}}'_{\Delta t} \times \overline{B}_{\Delta t}' + \frac{1}{\sigma} \overline{J}_{\Delta t} - \frac{\epsilon}{\mu ^2} \left| \operatorname{curl} \overline{B}_{\Delta t} \right|^2\operatorname{curl} \overline{B}_{\Delta t} \right] \cdot \operatorname{curl} b \nonumber \\
&- \epsilon \operatorname{curl}\left(\operatorname{curl} \overline{B}_{\Delta t}\right) \cdot \operatorname{curl}\left(\operatorname{curl} b\right)\ dxdt. \label{-38}
\end{align}

\subsection{Energy inequality on the \texorpdfstring{$\Delta t$-level}{}} \label{discreteenergyinequality}

In contrast to the incompressible setting in \cite{incompressiblecase} we have to combine the discrete induction equation \eqref{116} with the continuous Navier-Stokes equations \eqref{8}, \eqref{57} in a suitable way in order to derive an energy inequality at the $\Delta t$-level. We pick an arbitrary $t \in (0,T]$ and choose $k \in  \left\lbrace 1,...,\frac{T}{\Delta t} \right\rbrace$, $\xi \in [0,\Delta t)$ such that $t = k\Delta t - \xi$. For the magnetic part of the energy inequality we test the induction equation \eqref{116} by $\frac{1}{\mu}B_{\Delta t}^k$, which leads to
\begin{align}
& \frac{1}{2\mu}\partial_t B_{\Delta t, \|\cdot \|}(t) + \frac{1}{\mu} \int_\Omega \frac{1}{\sigma \mu} \left| \operatorname{curl}B_{\Delta t}^k \right|^2 + \frac{\epsilon}{\mu^2} \left| \operatorname{curl}B_{\Delta t}^k \right|^4 + \epsilon \left| \Delta B_{\Delta t}^k \right|^2\ dx \nonumber \\
=&\frac{1}{2\mu \Delta t} \left(\left\| B_{\Delta t}^k \right\|_{L^2(\Omega)}^2 - \left\| B_{\Delta t}^{k-1} \right\|_{L^2(\Omega)}^2\right) + \frac{1}{\mu} \int_\Omega \frac{1}{\sigma \mu} \left| \operatorname{curl}B_{\Delta t}^k \right|^2 + \frac{\epsilon}{\mu^2} \left| \operatorname{curl}B_{\Delta t}^k \right|^4 + \epsilon \left| \Delta B_{\Delta t}^k \right|^2\ dx \nonumber \\
\leq& \frac{1}{\mu}\int_\Omega \left( \tilde{u}_{\Delta t}^{k-1} \times B_{\Delta t}^{k-1} \right) \cdot \operatorname{curl} B_{\Delta t}^k + \frac{1}{\sigma}J_{\Delta t}^k \cdot \operatorname{curl}B_{\Delta t}^k\ dx. \label{-17}
\end{align}

Since corresponding estimates hold true also for all time indices $l=1,...,k-1$ we can integrate (discretely) over the interval $[0,t]$, which yields
\begin{align}
& \frac{1}{2\mu} B_{\Delta t,\| \cdot \|}(t) + \Delta t \sum
_{l=1}^{k-1} \frac{1}{\mu} \int_\Omega \frac{1}{\sigma \mu} \left| \operatorname{curl}B_{\Delta t}^l \right|^2 + \frac{\epsilon}{\mu^2} \left| \operatorname{curl}B_{\Delta t}^l \right|^4 + \epsilon \left| \Delta B_{\Delta t}^l \right|^2\ dx \nonumber \\
&+ \frac{\Delta t - \xi}{\mu} \int_\Omega \frac{1}{\sigma \mu} \left| \operatorname{curl}B_{\Delta t}^{k} \right|^2 + \frac{\epsilon}{\mu^2} \left| \operatorname{curl}B_{\Delta t}^{k} \right|^4 + \epsilon \left| \Delta B_{\Delta t}^{k} \right|^2\ dx \nonumber \\
\leq& \frac{1}{2\mu} B_{\Delta t,\| \cdot \|}(0) + \Delta t \sum_{l=1}^{k-1} \frac{1}{\mu}\int_\Omega \left( \tilde{u}_{\Delta t}^{l-1} \times B_{\Delta t}^{l-1} \right) \cdot \operatorname{curl} B_{\Delta t}^l + \frac{1}{\sigma}J_{\Delta t}^l \cdot \operatorname{curl}B_{\Delta t}^l\ dx \nonumber \\
&+ \frac{\Delta t - \xi}{\mu}\int_\Omega \left( \tilde{u}_{\Delta t}^{k-1} \times B_{\Delta t}^{k-1} \right) \cdot \operatorname{curl} B_{\Delta t}^{k} + \frac{1}{\sigma}J_{\Delta t}^{k} \cdot \operatorname{curl}B_{\Delta t}^{k}\ dx \nonumber \\
\leq& \frac{1}{2\mu} \int_\Omega \left| B_0 \right|^2\ dx + \Delta t \sum_{l=1}^{k-1} \left[ \frac{c}{\epsilon} \left\| B_{\Delta t}^{l-1} \right\|_{L^2(\Omega)}^2 + \frac{\epsilon}{8} \left\| \tilde{u}_{\Delta t}^{l-1} \right\|_{L^4(\Omega)}^4 + \frac{\epsilon}{8\mu^3} \left\| \operatorname{curl}B_{\Delta t}^l \right\|_{L^4(\Omega)}^4 + \frac{c}{\epsilon^\frac{1}{3}} \left\| J_{\Delta t}^l \right\|_{L^\frac{4}{3}(\Omega)}^\frac{4}{3} \right. \nonumber \\
 & \left. + \frac{\epsilon}{8\mu^3} \left\| \operatorname{curl}B_{\Delta t}^l \right\|_{L^4(\Omega)}^4 \right] + \left(\Delta t - \xi \right) \left[ \frac{c}{\epsilon} \left\| B_{\Delta t}^{k-1} \right\|_{L^2(\Omega)}^2 + \frac{\epsilon}{8} \left\| \tilde{u}_{\Delta t}^{k-1} \right\|_{L^4(\Omega)}^4 + \frac{\epsilon}{8\mu^3} \left\| \operatorname{curl}B_{\Delta t}^k \right\|_{L^4(\Omega)}^4 \right. \nonumber \\
 &\left. + \frac{c}{\epsilon^\frac{1}{3}} \left\| J_{\Delta t}^k \right\|_{L^\frac{4}{3}(\Omega)}^\frac{4}{3} + \frac{\epsilon}{8\mu^3} \left\| \operatorname{curl}B_{\Delta t}^k \right\|_{L^4(\Omega)}^4 \right] \label{-234}
\end{align}

under exploitation of Hölder's and Young's inequalities in the last estimate. On the right-hand side of \eqref{-234} we further estimate, due to the definition of $\tilde{u}_{\Delta t}^{l-1}$ in \eqref{-235} and Jensen's inequality,
\begin{align}
\left\| \tilde{u}_{\Delta t}^{l-1} \right\|_{L^4(\Omega)}^4 \leq \frac{1}{\Delta t} \int_\Omega \int_{(l-2)\Delta t}^{(l-1)\Delta t} \left| u_{\Delta t,l-1}(\tau) \right|^4\ d\tau dx. \nonumber
\end{align}

Moreover, a direct calculation yields
\begin{align}
\Delta t \sum_{l=1}^{k-1} \frac{c}{\epsilon} \left\| B_{\Delta t}^{l-1} \right\|_{L^2(\Omega)}^2 + \left(\Delta t - \xi \right) \frac{c}{\epsilon} \left\| B_{\Delta t}^{k-1} \right\|_{L^2(\Omega)}^2 \leq& \frac{c\Delta t}{\epsilon} \left( \frac{\left\| B_0 \right\|_{L^2(\Omega)}^2}{2} + \frac{\left\| B_{\Delta t}^{k-1} \right\|_{L^2(\Omega)}^2}{2} + \sum_{l=2}^{k-1} \left\| B_{\Delta t}^{l-1} \right\|_{L^2(\Omega)}^2 \right) \nonumber \\
=& \frac{c}{\epsilon} \int_0^{(k-1)\Delta t} B_{\Delta t,\| \cdot \|}(\tau)\ d\tau \leq \frac{c}{\epsilon} \int_0^{t} B_{\Delta t,\| \cdot \|}(\tau)\ d\tau. \label{-236}
\end{align}

Hence, absorbing the $\|\operatorname{curl} B_{\Delta t}^l\|_{L^4(\Omega)}^4$-terms in \eqref{-234} into the left-hand side and expressing the sums as integrals, we end up with
\begin{align}
&\frac{1}{2\mu} B_{\Delta t,\| \cdot \|}(t) + \frac{1}{\mu} \int_0^t \int_\Omega \frac{1}{\sigma \mu} \left| \operatorname{curl}\overline{B}_{\Delta t} \right|^2 + \frac{6\epsilon}{8\mu^2} \left| \operatorname{curl}\overline{B}_{\Delta t} \right|^4 + \epsilon \left| \Delta \overline{B}_{\Delta t} \right|^2\ dxd\tau \nonumber \\
\leq& \frac{1}{2\mu} \int_\Omega \left| B_0 \right|^2\ dx + \frac{c}{\epsilon^\frac{1}{3}} \int_0^T \int_\Omega \left| \overline{J}_{\Delta t} \right|^\frac{4}{3}\ dxd\tau + \frac{\epsilon}{8}\Delta t \int_\Omega \left| u_0 \right|^4\ dx + \frac{\epsilon}{8} \int_{0}^{(k-1)\Delta t} \int_\Omega \left| u_{\Delta t} \right|^4 \ d\tau dx \nonumber \\
&+ \frac{c}{\epsilon} \int_0^{t} B_{\Delta t,\| \cdot \|}(\tau)\ d\tau. \label{-20}
\end{align}

For the mechanical part of the energy inequality we test the continuity equation \eqref{-246} by $\frac{1}{2}|u_{\Delta t}|^2$, the momentum equation \eqref{-248} by $u_{\Delta t}$, add up the resulting equations and obtain
\begin{align}
&\int_\Omega \frac{1}{2} \rho_{\Delta t}(t)\left|u_{\Delta t}(t)\right|^2 + \frac{a}{\gamma -1}\rho_{\Delta t}^\gamma(t) + \frac{\alpha }{\beta - 1}\rho_{\Delta t}^\beta(t)\ dx + \int_{0}^{t}\int_\Omega 2\nu\left( \overline{\chi}_{\Delta t}' \right) |\mathbb{D} \left(u_{\Delta t}\right)|^2  \nonumber \\
&+ \lambda \left( \overline{\chi}_{\Delta t}' \right) \left| \operatorname{div}u_{\Delta t} \right|^2+ a\epsilon \gamma \rho_{\Delta t}^{\gamma - 2}\left| \nabla \rho_{\Delta t} \right|^2 + \alpha \epsilon \beta \rho_{\Delta t}^{\beta - 2}\left| \nabla \rho_{\Delta t} \right|^2 + \epsilon \left| u_{\Delta t} \right|^4 \ dxd\tau \nonumber \\
\leq& \int_\Omega \frac{1}{2} \rho_0\left|u_0\right|^2 + \frac{a}{\gamma -1}\rho_0^\gamma + \frac{\alpha }{\beta - 1}\rho_0^\beta\ dx + \int_{0}^t\int_\Omega \rho_{\Delta t}g \cdot u_{\Delta t} + \frac{1}{\mu}\left( \text{curl}\overline{B}_{\Delta t}' \times \overline{B}_{\Delta t}' \right) \cdot u_{\Delta t}\ dxd\tau \nonumber \\
\leq& \int_\Omega \frac{1}{2} \rho_0\left|u_0\right|^2 + \frac{a}{\gamma -1}\rho_0^\gamma + \frac{\alpha }{\beta - 1}\rho_0^\beta\ dx + \int_{0}^t\int_\Omega \rho_{\Delta t}g \cdot u_{\Delta t} \ dxd\tau \nonumber \\
& + \frac{\epsilon}{8} \int_0^t \int_\Omega \left| u_{\Delta t} \right|^4\ dxd\tau + \frac{c}{\epsilon} \int_0^t B_{\Delta t, \| \cdot \|} (\tau)d\tau + \int_0^t \int_\Omega \frac{\epsilon}{8\mu^3} \left| \operatorname{curl} \overline{B}_{\Delta t}' \right|^4\ dxd\tau, \label{-55}
\end{align}

where the last estimate uses Hölder's inequality, Young's inequality and the same inequality \eqref{-236} as in the estimate of the corresponding term in the induction equation. Adding the inequality \eqref{-55} to the inequality \eqref{-20} and absorbing multiple terms from the right-hand side into the left-hand side, we finally get the energy inequality
\begin{align}
&\int_\Omega \frac{1}{2} \rho_{\Delta t}(t)\left|u_{\Delta t}(t)\right|^2 + a \frac{\rho_{\Delta t}^\gamma(t)}{\gamma -1} + \frac{\alpha \rho_{\Delta t}^\beta(t)}{\beta - 1}\ dx + \frac{1}{2\mu} B_{\Delta t,\| \cdot \|}(t) + \int_{0}^{t}\int_\Omega 2\nu\left( \overline{\chi}_{\Delta t}' \right) |\mathbb{D} \left(u_{\Delta t}\right)|^2  \nonumber \\
&+ \lambda \left( \overline{\chi}_{\Delta t}' \right) \left| \operatorname{div}u_{\Delta t} \right|^2+ a\epsilon \gamma \rho_{\Delta t}^{\gamma - 2}\left| \nabla \rho_{\Delta t} \right|^2 + \alpha \epsilon \beta \rho_{\Delta t}^{\beta - 2}\left| \nabla \rho_{\Delta t} \right|^2 + \frac{3\epsilon}{4} \left| u_{\Delta t} \right|^4 \ dxdt \nonumber \\
&+ \frac{1}{\mu} \int_0^t \int_\Omega \frac{1}{\sigma \mu} \left| \operatorname{curl}\overline{B}_{\Delta t} \right|^2 + \frac{5\epsilon}{8\mu^2} \left| \operatorname{curl}\overline{B}_{\Delta t} \right|^4 + \epsilon \left| \Delta \overline{B}_{\Delta t} \right|^2\ dxd\tau \nonumber \\
\leq& \int_\Omega \frac{1}{2} \rho_0\left|u_0\right|^2 + \frac{a}{\gamma -1}\rho_0^\gamma + \frac{\alpha }{\beta - 1}\rho_0^\beta\ dx + \frac{1}{2\mu} \int_\Omega \left| B_0 \right|^2\ dx + \frac{\epsilon}{8}\Delta t \int_\Omega \left| u_0 \right|^4\ dx + \frac{c}{\epsilon^\frac{1}{3}} \int_0^T \int_\Omega \left| \overline{J}_{\Delta t} \right|^\frac{4}{3}\ dxd\tau \nonumber \\
&+\int_{0}^t\int_\Omega \rho_{\Delta t}g \cdot u_{\Delta t} \ dxd\tau + \int_0^t B_{\Delta t, \| \cdot \|} (\tau)d\tau \nonumber \\
\leq& c + \int_{0}^t\int_\Omega \rho_{\Delta t}g \cdot u_{\Delta t} \ dxd\tau + \int_0^t B_{\Delta t, \| \cdot \|} (\tau)d\tau \quad \quad \forall t \in [0,T], \label{-26}
\end{align}

where the constant $c >0$ is independent of $\Delta t$ and $t$. In particular, by use of the Gronwall Lemma and the estimates for the solution to the Neumann problem for the density, c.f.\ Lemma \ref{neumannproblem}, we find a constant $c>0$, independent of $\Delta t$, such that the following bounds hold true:
\begin{align}
\left\| u_{\Delta t} \right\|_{C([0,T];V_n)} + \left\| \rho_{\Delta t} \right\|_{C([0,T];C^2( \overline{\Omega} ))} + \left\| \partial_t \rho _{\Delta t} \right\|_{C(\overline{Q})} \leq& c, \label{-237} \\
\left\| B_{\Delta t} \right\|_{L^\infty (0,T;L^2(\Omega))} + \left\| B_{\Delta t} \right\|_{L^2 (0,T;H^{2,2}(\Omega))} + \left\| \operatorname{curl} B_{\Delta t} \right\|_{L^4 (Q)} \leq& c, \label{-238} \\
\left\| \overline{B}_{\Delta t} \right\|_{L^\infty (0,T;L^2(\Omega))} + \left\| \overline{B}_{\Delta t} \right\|_{L^2 (0,T;H^{2,2}(\Omega))} + \left\| \operatorname{curl} \overline{B}_{\Delta t} \right\|_{L^4 (Q)} \leq& c, \label{-239} \\
\left\| \overline{B}_{\Delta t}' \right\|_{L^\infty (0,T;L^2(\Omega))} + \left\| \overline{B}_{\Delta t}' \right\|_{L^2 (0,T;H^{2,2}(\Omega))} + \left\| \operatorname{curl} \overline{B}_{\Delta t}' \right\|_{L^4 (Q)} \leq& c. \label{-240}
\end{align}
The bounds for the magnetic induction in $L^\infty (0,T;L^2(\Omega))$ in \eqref{-238}--\eqref{-240} follow from the choice $t = k\Delta t$, $k=1,...,\frac{T}{\Delta t}$ in the energy inequality \eqref{-26}, for which it holds $B_{\Delta t,\| \cdot \|}(t) = \| B_{\Delta t}^k \|_{L^2(\Omega)}^2$. For a bound of the time derivative of $u_{\Delta t}$ we introduce the operator
\begin{align}
\mathcal{M}_{\rho_{\Delta t}(t)}: V_n \rightarrow V_n,\quad \quad \left\langle M_{\rho_{\Delta t}(t)}v^1,v^2 \right\rangle := \int_\Omega \rho_{\Delta t}(t) v^1 \cdot v^2\ dx \quad \quad \forall v^1,v^2 \in V_n \nonumber
\end{align}

and denote
\begin{align}
\mathcal{N}\left(\rho_{\Delta t}, u_{\Delta t} \right) :=& -\operatorname{div} \left( \rho_{\Delta t} u_{\Delta t} \otimes u_{\Delta t} \right) - a \nabla \rho_{\Delta t}^\gamma - \alpha \nabla \rho_{\Delta t}^\beta + \operatorname{div} \left( 2\nu \left( \chi_{\Delta t}^{k-1} \right)\mathbb{D}\left(u_{\Delta t}\right) \right) \nonumber \\
+& \nabla \left( \lambda \left( \chi_{\Delta t}^{k-1} \right)\operatorname{div}\left(u_{\Delta t}\right) \right) + \rho_{\Delta t} g + \frac{1}{\mu} \left( \text{curl}B_{\Delta t}^{k-1} \times B_{\Delta t}^{k-1} \right) - \epsilon \nabla u_{\Delta t} \nabla \rho_{\Delta t} - \epsilon \left|u_{\Delta t}\right|^2 u_{\Delta t}. \nonumber
\end{align}

Due to the uniform bound \eqref{-237} the solution $\rho_{\Delta t}$ to the Neumann problem \eqref{-246}, \eqref{-247} is bounded away from $0$ uniformly with respect to $\Delta t$, c.f.\ the estimate \eqref{-329} in Lemma \ref{neumannproblem} below. Consequently, $\mathcal{M}_{\rho_{\Delta t}}$ is invertible and from the momentum equation \eqref{-248} it is possible to derive the representation
\begin{align}
\partial_t u_{\Delta t}(t) =& \mathcal{M}_{\rho_{\Delta t}(t)}^{-1} \mathcal{M}_{\partial_t\rho_{\Delta t}(t)}\mathcal{M}_{\rho_{\Delta t}(t)}^{-1} \left[ P_n \left(\rho_0 u_0 \right) + \int_0^t P_n \left( \mathcal{N}\left(\rho_{\Delta t}, u_{\Delta t}\right) \right)dt \right] \nonumber \\
&+ \mathcal{M}_{\rho_{\Delta t}(t)}^{-1} \left[ P_n \left( \mathcal{N}\left(\rho_{\Delta t}, u_{\Delta t}\right) \right)(t) \right]. \nonumber
\end{align}

This together with the uniform bound of $\rho_{\Delta t}$ away from $0$ and the uniform bounds \eqref{-237}--\eqref{-240} leads to the estimate
\begin{align}
\left\| \partial_t u_{\Delta t} \right\|_{L^2(0,T;V_n)}^2 \leq c. \label{-251} 
\end{align}

For more details on the derivation of \eqref{-251} we refer to \cite[Section 7.7]{novotnystraskraba} and in particular \cite[Section 7.7.4.1]{novotnystraskraba}. The bounds \eqref{-237}--\eqref{-240} and \eqref{-251} and the Aubin-Lions Lemma imply the existence of functions
\begin{align}
0 \leq \rho &\in \left\lbrace \psi \in C\left([0,T];H^{1,2}(\Omega)\right) \bigcap L^2 \left(0,T; H^{2,2}(\Omega) \right):\ \partial_t \psi \in L^2 \left(Q \right),\ \left. \nabla \psi \cdot \operatorname{n} \right|_{\partial \Omega} = 0 \right\rbrace, \label{-354} \\
u &\in \left\lbrace \phi \in C\left([0,T];V_n \right):\ \partial_t \phi \in L^2\left(0,T;V_n \right) \right\rbrace, \nonumber \\
B &\in \left\lbrace b \in L^\infty \left(0,T;L^2(\Omega)\right) \bigcap L^2\left(0,T;H^{2,2}_{\operatorname{div}}(\Omega)\right):\ b \cdot \text{n}|_{\partial \Omega} = 0 \right\rbrace \label{-355}
\end{align}

such that, after the extraction of a subsequence, it holds
\begin{align}
\rho_{\Delta t} \rightharpoonup \rho \ \ \quad &\text{in } L^2\left(0,T;H^{2,2}(\Omega) \right), &\rho_{\Delta t}& \rightarrow \rho \ \quad \text{in } C\left([0,T];H^{1,2}(\Omega) \right), \label{-261} \\
\partial_t \rho_{\Delta t} \rightharpoonup \partial_t \rho \quad &\text{in } L^2\left(Q \right),\quad &\rho_{\Delta t}^\frac{\beta}{2} & \rightharpoonup \rho^\frac{\beta}{2} \quad \text{in } L^2\left(0,T; H^{1,2}(\Omega)\right), \label{-253} \\
\quad u_{\Delta t} \rightarrow u \ \ \quad &\text{in } C\left([0,T];V_n \right), \quad &\partial_t u_{\Delta t}& \rightharpoonup \partial_t u \quad \text{in } L^2\left(0,T;V_n \right), \label{-256} \\
B_{\Delta t},\ \overline{B}_{\Delta t},\ \overline{B}'_{\Delta t} \buildrel\ast\over\rightharpoonup B \quad &\text{in } L^\infty \left(0,T;L^2\left(\Omega \right)\right), \quad &B_{\Delta t},\ \overline{B}_{\Delta t},\ \overline{B}'_{\Delta t}& \rightharpoonup B \quad \ \  \text{in } L^2\left(0,T;H^{2,2}\left(\Omega \right)\right). \nonumber
\end{align}

Here, the fact that the weak limits of the different interpolants coincide follows from Lemma \ref{equalityofrothelimits}. Moreover the boundary conditions of the limit functions in \eqref{-354} and \eqref{-355} follow directly from the corresponding boundary conditions on the $\Delta t$-level, c.f.\ Proposition \ref{System on Delta t-level} and the definition of the space $Y^k(S_{\Delta t})$ in \eqref{-315}. Furthermore, the external force $\overline{J}_{\Delta t}$, discretized via \eqref{840}, converges to its original time-dependent counterpart,
\begin{align}
\overline{J}_{\Delta t} \rightarrow J \quad \text{in }L^q(Q) \quad \quad \quad \forall 1 \leq q < \infty. \nonumber
\end{align}

Finally, $X_{\Delta t}$, as the solution to the initial value problem \eqref{-254}, satisfies the conditions of Lemma \ref{convergence} which tells us that
\begin{align}
X_{\Delta t} \rightarrow X \quad &\text{in } C\left([0,T];C_{\text{loc}}\left( \mathbb{R}^3 \right) \right), \nonumber \\
\chi_{\Delta t} \rightarrow \chi := \textbf{db}_{S(\cdot)}(\cdot) \quad &\text{in } C\left([0,T];C_{\text{loc}}\left( \mathbb{R}^3 \right) \right), \label{-37}
\end{align}

where $S(t) := (X(t;O))^\delta$ and $X$ represents the solution to
\begin{align}
\frac{dX(t;x)}{dt} = R_\delta \left[ u \right] \left( t, X(t;x) \right),\quad \text{for } t \in [0,T], \quad \quad X(0;x) =  x,\quad \text{for } x \in \mathbb{R}^3. \nonumber
\end{align}

\subsection{Continuity equation}

Due to the convergences \eqref{-261}--\eqref{-256} of the density and the velocity we can pass to the limit in the continuity equation \eqref{-246} and infer that the limit functions $\rho$ and $u$ solve the initial value problem
\begin{align}
\partial_t \rho + \nabla \left( \rho u \right) =& \epsilon \Delta \rho \quad \text{a.e. in } Q,\quad \quad \rho(0) = \rho_0. \nonumber
\end{align}

\subsection{Induction equation} \label{rothelimitinduction}

We first show convergence of the quantity $\overline{\tilde{u}}_{\Delta t}'$ from the mixed term in the discrete induction equation \eqref{-38}. We fix an arbitrary point $(t,x) \in Q$ and, for each sufficiently small $\Delta t > 0$, we choose $k_{\Delta t} \in \lbrace 2,...,\frac{T}{\Delta t} \rbrace$ such that $t \in [(k_{\Delta t}-1)\Delta t, k_{\Delta t}\Delta t)$. It holds
\begin{align}
\left| \overline{\tilde{u}}_{\Delta t}'(t,x) - u(t,x) \right| =& \left| \frac{1}{\Delta t} \int_{(k_{\Delta t}-2)\Delta t}^{(k_{\Delta t}-1)\Delta t} u_{\Delta t,k_{\Delta t} - 1}(\tau,x) d\tau - u(t,x) \right| \nonumber \\
\leq& \sup_{\tau \in [(k_{\Delta t}-2)\Delta t, (k_{\Delta t}-1)\Delta t]} \left| u_{\Delta t}(\tau - \Delta t,x) - u(t,x) \right| \rightarrow 0 \label{-259}
\end{align}

due to the uniform convergence \eqref{-256} of $u_{\Delta t}$. Moreover, the uniform bound of $u_{\Delta t}$ in \eqref{-237} shows equiintegrability of $|\overline{\tilde{u}}_{\Delta t}' - u_{\Delta t}|^q$ for any $1 \leq q < \infty$. This together with the pointwise convergence \eqref{-259} gives us the conditions for the Vitali convergence theorem and we infer that
\begin{align}
\overline{\tilde{u}}'_{\Delta t} \rightarrow u \quad \text{in } L^q(Q) \quad \forall 1 \leq q < \infty. \label{-36}
\end{align}

Further, due to the uniform bounds \eqref{-239}, we find a function $z \in L^\frac{4}{3}(Q)$ such that, possibly after the extraction of a suitable subsequence, it holds
\begin{align}
\epsilon \left| \operatorname{curl}\overline{B}_{\Delta t} \right|^2 \operatorname{curl}\overline{B}_{\Delta t} \rightharpoonup \epsilon z \quad \text{in } L^\frac{4}{3} (Q). \label{-56}
\end{align}

We remark that, since the quantity $\epsilon z$ will vanish from the system as soon as we let $\epsilon$ tend to zero, there is no need to specify the form of the limit function $z$ in \eqref{-56}. Now we test \eqref{-38} by an arbitrary function $b \in Y(S)$. This is possible, since $b$ is curl-free in an open neighbourhood of $\overline{Q}^s(S)$ and so, by the uniform convergence \eqref{-37} of the signed distance function, it also satisfies
\begin{align}
\operatorname{curl} b(t) = 0 \quad \text{a.e. in } S_{\Delta t}(k\Delta t) \quad  \text{for a.a. } t \in [(k-1)\Delta t, k \Delta t],\ k = 1,...,\frac{T}{\Delta t} \label{-260}
\end{align}

for all sufficiently small $\Delta t>0$. Hence, for all such $t$, $k$ and $\Delta t$ it holds $b(t) \in W^k(S_{\Delta t})$ and thus $b$ is indeed an admissible test function in \eqref{-38}, c.f. \eqref{-41}. Using the strong convergence \eqref{-36} of $\overline{\tilde{u}}'_{\Delta t}$, we can pass to the limit with respect to $\Delta t \rightarrow 0$ in the resulting equation and obtain the induction equation in the limit,
\begin{align}
-\int_0^T \int_\Omega B \cdot \partial_t b\ dxdt =& \int_0^T\int_\Omega \left[ -\frac{1}{\sigma \mu} \operatorname{curl} B + u \times B + \frac{1}{\sigma} J - \frac{\epsilon}{\mu ^2} z \right] \cdot \operatorname{curl} b \nonumber \\
&- \epsilon \operatorname{curl}\left(\operatorname{curl} B\right) \cdot \operatorname{curl}\left(\operatorname{curl} b\right)\ dxdt \quad \quad \forall b \in Y(S). \nonumber
\end{align}

\subsection{Momentum equation}

In order to pass to the limit in the momentum equation, it remains to show convergence of the piecewise constant Lorentz force. This is achieved by the same arguments as in the incompressible case, c.f.\ \cite[Section 4.2]{incompressiblecase}: The uniform bounds \eqref{-239} and \eqref{-240} allow us to extract suitable subsequences and find functions $z_1,\ z_2 \in L^2(Q)$ with the properties that
\begin{align}
\operatorname{curl} \overline{B}'_{\Delta t} \times \overline{B}'_{\Delta t} \rightharpoonup z_1 \quad &\text{in } L^2(Q), \quad &\operatorname{curl} \overline{B}_{\Delta t} \times \overline{B}'_{\Delta t}& \rightharpoonup z_2 \quad \text{in } L^2(Q). \label{-54}
\end{align}

Our goal is to identify the limit functions $z_1$ and $z_2$ as
\begin{align}
z_1 = z_2 = \operatorname{curl} B \times B \quad \text{a.e. in } Q. \label{-49}
\end{align}

Since, according to \eqref{-318}, it holds
\begin{align}
\left| Q \right| = \left| Q^s(S) \bigcup Q^f(S) \right|, \nonumber
\end{align}

it is sufficient to do so in $Q^s(S)$ and $Q^f(S)$. Focussing at first on the fluid part of the domain, we consider any $c,d \in [0,T]$ and any ball $U \subset \Omega$ such that $\overline{I \times U} := \overline{(c,d)\times U} \subset Q^f(S)$. In particular, any function $b \in L^4(I;H_0^{2,2}(U))$, extended by $0$ outside of $\overline{I \times U}$, is curl-free in an open neighbourhood of $\overline{Q}^s(S)$. Hence, the uniform convergence \eqref{-37} of the signed distance function implies that, for sufficiently small $\Delta t > 0$, any such $b$ also satisfies the curl-free condition \eqref{-260} in $S_{\Delta t}(k\Delta t)$, $k=1,...,\frac{T}{\Delta t}$, and therefore constitutes an admissible test function in the induction equation \eqref{-38}. Using it as such we deduce the uniform bound
\begin{align}
\left\| \frac{\overline{B}'_{\Delta t}(\cdot) - \overline{B}'_{\Delta t}(\cdot -\Delta t)}{\Delta t} \right\|_{L^\frac{4}{3}\left(a + \Delta t, d;H^{-2,2}(U)\right)} \leq \left\| \frac{\overline{B}_{\Delta t}(\cdot) - \overline{B}_{\Delta t}(\cdot -\Delta t)}{\Delta t} \right\|_{L^\frac{4}{3}(a,d;H^{-2,2}(U))}\leq c.\nonumber
\end{align}

Since, by \eqref{-240}, $\overline{B}'_{\Delta t}$ is bounded uniformly in $L^2(Q)$, we can thus apply the discrete Aubin-Lions Lemma \cite[Theorem 1]{dreherjungel} and conclude that
\begin{align}
\overline{B}'_{\Delta t} \rightarrow B \quad \text{in } L^2\left(I;H^{-1,2}(U) \right). \nonumber
\end{align}

This strong convergence implies that
\begin{align}
z_1 = z_2 = \operatorname{curl} B \times B \quad \text{a.e. in } I \times U \ \text{and hence in } Q^f(S). \label{-346}
\end{align}

In order to show the equation \eqref{-49} also in the solid region, we consider another arbitrary pair of an interval $I \subset (0,T)$ and a ball $U \subset \Omega$, this time satisfying $\overline{I \times U} \subset Q^s(S)$. From the fact that $\operatorname{curl} B_{\Delta t}^k = 0$ in $S_{\Delta t}^k \bigcap \Omega$ for $k=1,...,\frac{T}{\Delta t}$, and the uniform convergence \eqref{-37} of the signed distance function it follows that, for any sufficiently small $\Delta t > 0$,
\begin{align}
\operatorname{curl} \overline{B}'_{\Delta t} = \operatorname{curl} \overline{B}_{\Delta t} = 0 \quad \text{a.e. in } I \times U. \nonumber
\end{align}

Thus, letting $\Delta t$ tend to zero, we infer that
\begin{align}
z_1 = z_2 = 0 = \operatorname{curl} B = \operatorname{curl} B \times B \quad  \text{a.e. in } I \times U\ \text{and hence in } Q^s(S). \label{-47}
\end{align}

In combination with the corresponding equality \eqref{-346} in $Q^f(S)$, we have therefore shown the desired identification \eqref{-49} of $z_1$ and $z_2$. We remark that \eqref{-47} moreover shows that the solid region $Q^s(S)$ in the limit is again insulating. Next we exploit the uniform convergence \eqref{-37} of the signed distance function together with the definition \eqref{-52} of the variable viscosity coefficients as smooth functions of $\overline{\chi}_{\Delta t}'$ to infer that
\begin{align}
\nu \left( \overline{\chi}'_{\Delta t} \right) \rightarrow \nu \left( \chi \right) \quad \text{in } C\left([0,T];C_\text{loc}\left( \mathbb{R}^3 \right)\right), \quad \quad \lambda \left( \overline{\chi}'_{\Delta t} \right) \rightarrow \lambda \left( \chi \right) \quad \text{in } C\left([0,T];C_\text{loc}\left( \mathbb{R}^3 \right)\right). \label{-58}
\end{align}

This, in combination with the convergence of the Lorentz force, c.f.\ \eqref{-54}, \eqref{-49} and the convergences \eqref{-261}--\eqref{-256} of $\rho_{\Delta t}$ and $u_{\Delta t}$, allows us to pass to the limit in the momentum equation \eqref{-248} and infer the equation
\begin{align}
&\int_{0}^{T} \int_\Omega \partial_t \left( \rho u \right) \cdot \phi\ dxdt \nonumber \\
=& \int _{0}^{T} \int _\Omega \left( \rho u \otimes u \right): \mathbb{D} (\phi) + \left( a\rho^\gamma + \alpha \rho^\beta \right) \operatorname{div} \phi - 2\nu \left( \chi \right)\mathbb{D}\left(u \right):\mathbb{D} (\phi) - \lambda \left( \chi \right)\operatorname{div}\left(u\right)\operatorname{div} \phi \nonumber \\
&+ \rho g\cdot \phi + \frac{1}{\mu}\left( \text{curl}B \times B \right) \cdot \phi - \epsilon \left| u \right|^2u \cdot \phi - \epsilon \left( \nabla u \nabla \rho \right) \cdot \phi\ dxdt \quad \forall \phi \in C\left([0,T]; V_n \right). \nonumber
\end{align}

\subsection{Limit passage in the energy inequality}

We choose an arbitrary $t \in (0,T]$ and $k \in \lbrace 1,...,\frac{T}{\Delta t} \rbrace$, such that $t = k\Delta t - \xi$ for some $\xi \in [0,\Delta t)$. Subsequently, we sum the inequality \eqref{-17} over $l=1,...,k$ and add the first inequality in \eqref{-55} to obtain
\begin{align}
&\int_\Omega \frac{1}{2}\rho_{\Delta t}(t)\left|u_{\Delta t}(t)\right|^2 + \frac{a}{\gamma -1}\rho_{\Delta t}^\gamma(t) + \frac{\alpha}{\beta - 1} \rho_{\Delta t}^\beta(t) + \frac{1}{2\mu} \left| \overline{B}_{\Delta t}(t) \right|^2\ dx \nonumber \\
&+ \int_{0}^{t}\int_\Omega 2\nu\left( \overline{\chi}_{\Delta t}' \right) |\mathbb{D} \left(u_{\Delta t}\right)|^2 + \lambda \left( \overline{\chi}_{\Delta t}' \right) \left| \operatorname{div}u_{\Delta t} \right|^2 + a\epsilon \gamma \rho_{\Delta t}^{\gamma - 2}\left| \nabla \rho_{\Delta t} \right|^2 + \alpha \epsilon \beta \rho_{\Delta t}^{\beta - 2}\left| \nabla \rho_{\Delta t} \right|^2 + \epsilon \left| u_{\Delta t} \right|^4 \nonumber \\
&+ \frac{1}{\sigma \mu^2} \left| \operatorname{curl} \overline{B}_{\Delta t} \right|^2 + \frac{\epsilon}{\mu^3} \left| \operatorname{curl} \overline{B}_{\Delta t} \right|^4 + \frac{\epsilon}{\mu} \left| \Delta \overline{B}_{\Delta t} \right|^2\ dxd\tau \nonumber \\
\leq& \int_\Omega \frac{1}{2}\rho_0 \left|u_0\right|^2 + \frac{a}{\gamma -1}\rho_0^\gamma + \frac{\alpha}{\beta - 1} \rho_0^\beta + \frac{1}{2\mu} \left| B_0 \right|^2\ dx + \int_{0}^t\int_\Omega \rho_{\Delta t}g \cdot u_{\Delta t} + \frac{1}{\mu}\left( \text{curl}\overline{B}_{\Delta t}' \times \overline{B}_{\Delta t}' \right) \cdot u_{\Delta t} \nonumber \\
&+ \frac{1}{\mu} \left( \overline{\tilde{u}}_{\Delta t}' \times \overline{B}_{\Delta t}' \right) \cdot \operatorname{curl} \overline{B}_{\Delta t} + \frac{1}{\sigma \mu} \overline{J}_{\Delta t} \cdot \operatorname{curl} \overline{B}_{\Delta t} \ dxd\tau + c \left( \Delta t \right)^\frac{1}{2}. \label{-65}
\end{align}

Here, on the left-hand side we drop the term $\rho_{\Delta t}^{\gamma - 2}\left| \nabla \rho_{\Delta t} \right|^2$ and pass to the limit by exploiting, in particular, the weak lower semicontinuity of norms and the strong convergence \eqref{-58} of the variable viscosity coefficients. Moreover, on the right-hand side of \eqref{-65} we can carry out the limit passage by exploiting the convergence \eqref{-54}, \eqref{-49} of the Lorentz force. Altogether we obtain
\begin{align}
&\int_\Omega \frac{1}{2}\rho(t)\left|u(t)\right|^2 + \frac{a}{\gamma -1}\rho^\gamma(t) + \frac{\alpha}{\beta - 1} \rho^\beta(t) + \frac{1}{2\mu} \left| B(t) \right|^2\ dx + \int_{0}^{t}\int_\Omega 2\nu\left( \chi \right) |\mathbb{D} \left(u\right)|^2 + \lambda \left( \chi \right) \left| \operatorname{div}u \right|^2 \nonumber \\
&+ \alpha \epsilon \beta \rho^{\beta - 2}\left| \nabla \rho \right|^2 + \epsilon \left| u \right|^4 + \frac{1}{\sigma \mu^2} \left| \operatorname{curl} B \right|^2 + \frac{\epsilon}{\mu} \left| \Delta B \right|^2 + \frac{\epsilon}{\mu^3} \left| z \right|^4 \ dxd\tau \nonumber \\
\leq& \int_\Omega \frac{1}{2}\rho_0 \left|u_0\right|^2 + \frac{a}{\gamma -1}\rho_0^\gamma + \frac{\alpha}{\beta - 1} \rho_0^\beta + \frac{1}{2\mu} \left| B_0 \right|^2\ dx + \int_{0}^t\int_\Omega \rho g \cdot u + \frac{1}{\sigma \mu} J \cdot \operatorname{curl} B \ dxd\tau \quad \text{for } t \in [0,T]. \nonumber
\end{align}

Hence we have proved

\begin{proposition}
\label{System on Galerkin-level}
Let $n \in \mathbb{N}$, $\eta, \epsilon, \alpha > 0$, $\beta > \max\lbrace 4, \gamma \rbrace$ sufficiently large and let the assumptions of Theorem \ref{mainresult} be satisfied. Furthermore, let
\begin{align}
\rho _0 \in& C^{2,\frac{1}{2}}\left(\overline{\Omega}\right),\quad \quad (\rho u)_0 \in C^2\left(\overline{\Omega}\right), \quad \quad u_0 := P_n \left(\frac{(\rho u)_0}{\rho_0}\right) \in V_n, \quad \quad B_0 \in H^{2,2}(\Omega), \nonumber \\
&0 < \alpha \leq \rho_0 \leq \alpha^{-\frac{1}{2\beta}},\quad \quad \left. \nabla \rho_0 \cdot \operatorname{n} \right|_{\partial \Omega} = 0, \quad \quad \operatorname{div} B_0 = 0,\quad \quad \left. B_0 \cdot \operatorname{n}\right|_{\partial \Omega} = 0. \nonumber
\end{align}

Under these conditions there exist functions $X_n: \mathbb{R}^3 \rightarrow \mathbb{R}^3$, $z_n \in L^\frac{4}{3} ((0,T)\times \Omega)$ and 
\begin{align}
0 \leq \rho_n &\in \left\lbrace \psi \in C\left([0,T];H^{1,2}(\Omega)\right) \bigcap L^2 \left(0,T; H^{2,2}(\Omega) \right):\ \partial_t \psi \in L^2 \left(Q \right),\ \left. \nabla \psi \cdot \operatorname{n} \right|_{\partial \Omega} = 0 \right\rbrace, \label{-359} \\
u_n &\in \left\lbrace \phi \in C\left([0,T];V_n \right):\ \partial_t \phi \in L^2\left(0,T;V_n \right) \right\rbrace, \nonumber \\
B_n &\in \left\lbrace b \in L^\infty \left(0,T;L^2(\Omega)\right) \bigcap L^2\left(0,T;H^{2,2}_{\operatorname{div}}(\Omega)\right):\ \operatorname{curl}b = 0 \ \text{in } Q^s(S_n),\ b \cdot \operatorname{n}|_{\partial \Omega} = 0 \right\rbrace \label{-360}
\end{align}

for $S_n = S_n(\cdot) = (X_n(\cdot;O))^\delta$, which satisfy
\begin{align}
\frac{dX_n(t;x)}{dt} = R_\delta \left[ u_n \right] \left( t, X_n(t;x) \right),\quad &\text{for } t \in [0,T], \quad \quad X_n(0;x) =  x,\quad \text{for } x \in \mathbb{R}^3, \nonumber \\
\partial_t \rho_n + \operatorname{div} \left( \rho_n u_n \right) &= \epsilon \Delta \rho_n \quad \text{a.e. in } Q \label{-267}
\end{align}

and
\begin{align}
\int_0^T \int_\Omega \partial_t \left( \rho_n u_n \right) \cdot \phi \ dxdt =& \int _0^T \int _\Omega \left( \rho_n u_n \otimes u_n \right): \mathbb{D} (\phi) + \left( a\rho_n^\gamma  + \alpha \rho_n^\beta \right) \operatorname{div} \phi - 2\nu \left( \chi_n \right) \mathbb{D}(u_n):\mathbb{D} (\phi) \nonumber \\
&- \lambda \left( \chi_n \right) \operatorname{div} u_n \operatorname{div} \phi + \rho_n g\cdot \phi + \frac{1}{\mu}\left( \operatorname{curl}B_n \times B_n \right) \cdot \phi - \epsilon \left| u_n \right|^2u_n \cdot \phi \nonumber \\
&- \epsilon \left( \nabla u_n \nabla \rho_n \right) \cdot \phi \ dxdt, \label{-268} \\
- \int _0^T \int_\Omega B_n \cdot \partial_t b\ dxdt =& \int _0^T \int _\Omega \left[ -\frac{1}{\sigma \mu} \operatorname{curl} B_n + u_n \times B_n + \frac{1}{\sigma} J - \frac{\epsilon}{\mu ^2} z_n \right] \cdot \operatorname{curl} b \nonumber \\
&- \epsilon \operatorname{curl}\left( \operatorname{curl} B_n \right): \operatorname{curl}\left( \operatorname{curl} b\right)\ dxdt, \label{-269}
\end{align}

where $\chi_{n}(t,x) := \textbf{db}_{S_{n}(t)}(x)$, for any $\phi \in C([0,T]; V_n )$ and any $b\in Y(S_n)$. Further, these functions satisfy the initial conditions
\begin{align}
\rho_n(0) = \rho_0,\quad \quad u_n(0) = u_0,\quad \quad B_n(0) = B_0, \nonumber
\end{align}

of which the latter identity can be understood in the sense of \eqref{-341}, as well as the energy inequality
\begin{align}
&\int_\Omega \frac{1}{2} \rho_n (t)\left|u_n (t)\right|^2 + \frac{a}{\gamma - 1} \rho_n ^\gamma(t) + \frac{\alpha}{\beta - 1}\rho_n ^\beta(t) + \frac{1}{2\mu} \left| B_n(t) \right|^2\ dx + \int_0^t \int_\Omega 2\nu \left(\chi_n \right) |\mathbb{D} \left(u_n \right)|^2 \nonumber \\
&+ \lambda \left(\chi_n \right) \left| \operatorname{div}u_n \right|^2 + \alpha \epsilon \beta \rho_n^{\beta - 2}\left| \nabla \rho_n \right|^2 + \epsilon \left| u_n \right|^4 + \frac{1}{\sigma \mu^2} \left| \operatorname{curl} B_n \right|^2 + \frac{\epsilon}{\mu} \left| \Delta B_n \right|^2 + \frac{\epsilon}{\mu^3} \left| z_n \right|^\frac{4}{3} \ dxd\tau \nonumber \\
\leq& \int_\Omega  \frac{1}{2} \rho_0\left|u_0\right|^2 + \frac{a}{\gamma - 1} \rho_0 ^\gamma + \frac{\alpha}{\beta - 1}\rho_0 ^\beta + \frac{1}{2\mu} \left| B_0 \right|^2\ dx +\int_{0}^t\int_\Omega \rho_n g \cdot u_n + \frac{1}{\sigma \mu} J \cdot \operatorname{curl} B_n \ dxd\tau \label{-271}
\end{align}

for almost all $t \in [0,T]$.

\end{proposition}

{\centering \section{Limit passage in the Galerkin method} \label{galerkinlimit} \par }

Next, we pass to the limit with respect to $n \rightarrow \infty$. Using Lebesgue interpolation, we infer from the energy inequality \eqref{-271} the existence of a constant $c(\epsilon, \alpha) > 0$, independent of $n$, such that
\begin{align}
\left\| \rho_n \right\|_{L^{\frac{5}{3}\beta}(Q)} \leq \left\| \rho_n \right\|^\frac{2}{5}_{L^\infty \left(0,T;L^\beta (\Omega)\right)} \left\| \rho_n \right\|^\frac{3}{5}_{L^\beta \left(0,T;L^{3\beta} (\Omega)\right)} \leq c(\epsilon, \alpha). \label{-69}
\end{align}

Moreover, from the classical $L^p$-$L^q$ regularity results for parabolic equations, c.f.\ \cite[Lemma 7.37, Lemma 7.38, Section 7.8.2]{novotnystraskraba}, we infer that $\rho_n$ as the solution to the regularized continuity equation \eqref{-267} satisfies the estimates
\begin{align}
\epsilon \left\| \nabla \rho_n \right\|_{L^{r}(Q)} + \epsilon \left\| \partial_t \rho_n \right\|_{L^{\tilde{r}}(Q)} + \epsilon^2 \left\| \Delta \rho_n \right\|_{L^{\tilde{r}}(Q)} \leq c \label{-66}
\end{align}

for
\begin{align}
r := \frac{10\beta - 6}{3\beta + 3} > 2,\quad \tilde{r} := \frac{5\beta - 3}{4\beta} > 1 \quad \quad \forall \beta > 6 \label{-347} 
\end{align}

and a constant $c>0$ independent of $n$, $\eta$ and $\epsilon$. The uniform bounds \eqref{-69},\eqref{-66} and the energy inequality \eqref{-271}, together with the Aubin-Lions Lemma, allow us to extract suitable subsequences and find functions $z \in L^\frac{4}{3}(Q)$ and
\begin{align}
u \in& L^2(0,T;H_0^{1,2}(\Omega)), \label{-356} \\
0 \leq \rho \in& \bigg\lbrace \psi \in L^\infty\left(0,T;L^\beta(\Omega)\right) \bigcap L^{r} \left(0,T; W^{1,r}(\Omega) \right) \bigcap L^{\tilde{r}}\left(0,T; W^{2,\tilde{r}}(\Omega) \right): \nonumber \\
&\partial_t \psi \in L^{\tilde{r}} \left((0,T) \times \Omega \right),\ \left. \nabla \psi \cdot \operatorname{n} \right|_{\partial \Omega} = 0 \bigg\rbrace \subset C\left([0,T];L^2(\Omega) \right), \label{-357} \\
B \in& \bigg\lbrace b \in L^\infty \left(0,T;L^2(\Omega)\right) \bigcap L^2\left(0,T;H^{2,2}_{\operatorname{div}}(\Omega)\right):\ b \cdot \text{n}|_{\partial \Omega} = 0 \bigg\rbrace \label{-358}
\end{align}

with the properties that
\begin{align}
\rho_n \rightarrow \rho \quad &\text{in } L^\beta \left(Q \right),\quad \quad &\rho_n& \rightarrow \rho \quad \ \ \text{in } L^2\left(0,T;H^{1,2}(\Omega) \right) \label{-274} \\
\rho_n \rightharpoonup \rho \quad &\text{in } L^{\tilde{r}}\left(0,T;W^{2,\tilde{r}}(\Omega) \right),\quad \quad &\partial_t \rho_n& \rightharpoonup \partial_t \rho \quad \text{in } L^{\tilde{r}}\left(Q \right), \label{-275} \\
u_n \rightharpoonup u \quad &\text{in } L^2 (0,T;H^{1,2} (\Omega)),\quad \quad &B_n& \buildrel\ast\over\rightharpoonup B \quad \ \ \text{in } L^\infty (0,T;L^{2} (\Omega)), \label{-72} \\
B_n \rightharpoonup B \quad &\text{in } L^2 (0,T;H^{2,2} (\Omega)),\quad \quad &z_n& \rightharpoonup z \quad \quad \text{in } L^\frac{4}{3}\left(Q \right). \nonumber
\end{align}

The boundary conditions of the limit functions in \eqref{-356}--\eqref{-358} follow directly from the corresponding boundary conditions \eqref{-359} and \eqref{-360} of $\nabla \rho_n$ and $B_n$ and the fact that $u_n \in V_n$ vanishes on $\partial \Omega$ for all $n \in \mathbb{N}$. Finally, the initial value problem \eqref{-254}, solved by $X_n$, yields that the conditions of Lemma \ref{convergence} are satisfied. Hence
\begin{align}
X_n \rightarrow X \quad &\text{in } C\left([0,T];C_{\text{loc}}\left( \mathbb{R}^3 \right) \right), \nonumber \\
\chi_n \rightarrow \chi := \textbf{db}_{S(\cdot)}(\cdot) \quad &\text{in } C\left([0,T];C_{\text{loc}}\left( \mathbb{R}^3 \right) \right), \label{-278}
\end{align}

where $S(t) := (X(t;O))^\delta$ and $X$ denotes the unique solution to the initial value problem
\begin{align}
\frac{dX(t;x)}{dt} = R_\delta \left[ u \right] \left( t, X(t;x) \right),\quad \text{for } t \in [0,T], \quad \quad X(0;x) =  x,\quad \text{for } x \in \mathbb{R}^3. \nonumber
\end{align}

\subsection{Continuity equation}

The convergences \eqref{-274}--\eqref{-72} of $\rho_n$ and $u_n$ allow us to pass to the limit in the continuity equation \eqref{-267}. Consequently, the limit functions $\rho$ and $u$ satisfy the continuity equation
\begin{align}
\partial_t \rho + \operatorname{div} \left( \rho u \right) =& \epsilon \Delta \rho \quad \text{a.e. in } Q,\quad \rho(0)=\rho_0. \nonumber
\end{align}

\subsection{Induction equation} \label{inductionequationetalevel}

At this stage - as well as in the later sections - the limit passage in the induction equation does not differ from the incompressible case, c.f.\ \cite{incompressiblecase}. For the convenience of the reader, we present the arguments here: We begin by making sure that test functions $b \in Y(S,X)$ for the limit equation are also admissible on the $n$-level. To this end we fix some arbitrary $\omega > 0$. Then the uniform convergence \eqref{-278} of the signed distance function implies the existence of $N = N (\omega) > 0$, such that for any function $b \in Y(S)$ with $\operatorname{curl} b = 0$ in an $\omega$-neighbourhood of $\overline{Q}^s(S)$ it also holds
\begin{align}
\operatorname{curl}b = 0 \quad \text{in an } \frac{\omega}{2}\text{-neighbourhood of } \overline{Q}^s(S_n) \quad \text{and thus}\quad b \in Y\left(S_n \right) \quad \forall n \geq N. \nonumber
\end{align}

In particular, for any interval $I \subset (0,T)$ and any ball $U \subset \Omega$, such that $\overline{I \times U} \subset Q^f(S)$ there exists $N = N (I \times U) > 0$ such that
\begin{align}
b \in Y\left( S_n \right) \quad \quad \forall n \geq N,\ b \in \mathcal{D}(I \times U), \label{-287}
\end{align}

where $b$ has been extended by $0$ outside of $I \times U$. Next, an interpolation between $L^\infty (0,T;L^2(\Omega))$ and $L^2(0,T;L^6(\Omega))$ provides a uniform bound of $B_n$ in $L^3(Q)$. In combination with the bounds of $u_n$ and $\operatorname{curl}B_n$ in $L^2(Q)$ this yields the existence of functions $z_3, z_4 \in L^\frac{6}{5}(Q)$ such that, for a suitable subsequence, it holds
\begin{align}
\operatorname{curl} B_n \times B_n \rightharpoonup z_3 \quad \text{in } L^\frac{6}{5}\left(Q \right), \quad \quad u_n \times B_n \rightharpoonup z_4 \quad \text{in } L^\frac{6}{5}\left(Q \right). \label{-75}
\end{align}

With the aim of identifying $z_3$ and $z_4$ we pick an arbitrary interval $I \subset (0,T)$ and an arbitrary ball $U \subset \Omega$ with the property $\overline{I \times U} \subset Q^f(S)$. From \eqref{-287} we know that for any sufficiently large $n \in \mathbb{N}$ the induction equation \eqref{-269} may be tested by all functions of the form $\psi b$, where $\psi \in \mathcal{D}(I)$, $b \in \mathcal{D}(U)$. Doing so, we obtain the dual estimate
\begin{align}
\left\| \partial_t \int_U B_n \cdot b dx \right\|_{L^\frac{4}{3}(I)} \leq c \nonumber
\end{align}

with a constant $c>0$ depending on $b$ but not on $n$. This, together with the Arzelá-Ascoli theorem, yields
\begin{align}
B_n \rightarrow B \quad \text{in } C_\text{weak}\left( \overline{I};L^2(U) \right) \quad \text{and hence in } L^2 \left( I;H^{-1,2}(U) \right). \label{-282}
\end{align}

Combining this with the weak convergences of $u_n$ and $B_n$ in $L^2(0,T;H^{1,2}(\Omega))$ we realize that
\begin{align}
z_3 = \operatorname{curl}B \times B ,\quad z_4 = u \times B \quad \quad \text{a.e. in } Q^f(S). \nonumber
\end{align}

Due to the uniform convergence \eqref{-278} of the signed distance function and the fact that $\operatorname{curl}B_n = 0$ in $Q^s(S_n)$, the former one of these identities also holds true in $Q^s(S)$,
\begin{align}
z_3 = 0 = \operatorname{curl}B = \operatorname{curl}B \times B \quad \text{a.e. in } Q^s(S). \nonumber
\end{align}

Since moreover test functions $b \in Y(S)$ are curl-free in $Q^s(S)$, we end up with
\begin{align}
z_3 = \operatorname{curl}B \times B ,\quad z_4 \cdot \operatorname{curl}b = u \times B \cdot \operatorname{curl}b \quad \quad \text{a.e. in } Q \label{-283}
\end{align}

for any $b \in Y(S)$. The relations \eqref{-75} and \eqref{-283} allow us to pass to the limit in the mixed term of the induction equations and hence, letting $n$ tend to infinity in \eqref{-269}, we see that
\begin{align}
-\int_0^T \int_\Omega B \cdot \partial_t b\ dxdt =& \int_0^T\int_\Omega \left[ -\frac{1}{\sigma \mu} \operatorname{curl} B + u \times B + \frac{1}{\sigma}J - \frac{\epsilon}{\mu ^2} z \right] \cdot \operatorname{curl} b \nonumber \\
&- \epsilon \operatorname{curl}\left(\operatorname{curl} B\right) \cdot \operatorname{curl}\left(\operatorname{curl} b\right)\ dxdt \nonumber
\end{align}

for any $b \in Y(S)$. Finally, for any $b \in \mathcal{D}(\Omega)$ with $\operatorname{curl}b = 0$ in a neighbourhood of $\overline{S}_0$, we can argue similarly as in the derivation of the $C_\text{weak}$-convergence \eqref{-282} of $B_n$ to deduce that
\begin{align}
\int_\Omega B_n(\cdot,x) \cdot b(x)\ dx \rightarrow \int_\Omega B(\cdot,x) \cdot b(x)\ dx \quad \text{in } C\left(\left[0,t_0\right] \right), \nonumber
\end{align}

for some small $t_0 = t_0(b)>0$. This yields the initial condition $B(0)=B_0$ in the sense of \eqref{-341}.

\subsection{Momentum equation} \label{dwfw3fef}

By the same methods as for the (purely mechanical) compressible Navier-Stokes system, c.f.\ \cite[Section 7.8.2]{novotnystraskraba}, we derive strong convergence of the momentum function in the Galerkin limit: Recalling that $P_n$ denotes the orthogonal projection of $L^2(\Omega)$ onto $V_n$, we test the momentum equation \eqref{-268} --- after a density argument --- by $P_n(\phi)$ for an arbitrary function $\phi \in L^\frac{5\beta - 3}{\beta-3}(0,T;H_0^{2,2}(\Omega))$. Since in particular
\begin{align}
\left| \int_0^T \int_\Omega \epsilon \left| u_n \right|^2u_n \cdot P_n(\phi)\ dxdt \right| \leq \epsilon \left\| \left| u_n \right|^3 \right\|_{L^\frac{4}{3}(Q)} \left\| P_n(\phi) \right\|_{L^4(Q)} \leq \epsilon \left\| u_n \right\|^3_{L^4(Q)} \left\| \phi \right\|_{L^{\frac{5\beta - 3}{\beta-3}}(0,T;H^{2,2}(\Omega))}, \nonumber
\end{align}

this results in the dual estimate
\begin{align}
\left\| \partial_t P_n \left( \rho_n u_n \right) \right\|_{L^{\tilde{r}}\left(0,T; H^{-2,2}(\Omega)\right)} = \left\| \partial_t P_n \left( \rho_n u_n \right) \right\|_{L^{\frac{5\beta - 3}{4\beta}}\left(0,T; H^{-2,2}(\Omega)\right)} \leq c. \nonumber
\end{align}

Consequently, from the Aubin-Lions Lemma, we conclude that
\begin{align}
P_n(\rho_n u_n) \rightarrow \rho u \quad \text{in } L^2\left( 0,T;H^{-1,2}(\Omega) \right) \quad \text{and hence} \quad \rho_n u_n \rightarrow \rho u \quad \text{in } L^2\left( 0,T;H^{-1,2}(\Omega) \right). \nonumber
\end{align}

Therefore, as $u_n$ converges weakly in $L^2(0,T;H^{1,2}(\Omega))$, we infer that, for example,
\begin{align}
\rho_n u_n \otimes u_n \rightharpoonup \rho u \otimes u \quad \text{in } L^2\left( 0,T; L^\frac{6\beta}{4\beta + 3} (\Omega) \right). \nonumber
\end{align}

Moreover, the bound of $u_n$ in $L^4(Q)$ implies the existence of some $\tilde{z} \in L^\frac{4}{3}(Q)$ such that, for a chosen subsequence, it holds
\begin{align}
\epsilon \left|u_n\right|^2 u_n \rightharpoonup \epsilon \tilde{z} \quad \text{in } L^\frac{4}{3}(Q). \nonumber
\end{align}

Using further the uniform convergence \eqref{-278} of the signed distance function for passing to the limit in the variable viscosity coefficients and the relations \eqref{-75}, \eqref{-283} which allow us to pass to the limit in the Lorentz force, we can now let $n$ tend to infinity in the momentum equation \eqref{-268} and infer that
\begin{align}
-\int_0^T \int_\Omega \rho u \cdot \partial_t \phi \ dxdt =& \int _0^T \int _\Omega \left( \rho u \otimes u \right): \mathbb{D} (\phi) + \left( a\rho^\gamma + \alpha \rho^\beta \right) \operatorname{div} \phi - 2\nu \left( \chi \right) \mathbb{D}(u):\mathbb{D} (\phi) \nonumber \\
&- \lambda \left( \chi \right) \operatorname{div} u \operatorname{div} \phi + \rho g\cdot \phi + \frac{1}{\mu}\left( \operatorname{curl}B \times B \right) \cdot \phi - \epsilon \tilde{z} \cdot \phi - \epsilon \left( \nabla u \nabla \rho \right) \cdot \phi \ dxdt \label{-286}
\end{align}

for any $\phi \in C_0^1 ([0,T]; V_N )$ with fixed $N \in \mathbb{N}$. Since $\bigcup_{n=1}^\infty V_n$ is dense in $H_0^{1,2}(\Omega)$, we finally conclude that \eqref{-286} also holds true for any $\phi \in \mathcal{D}(Q)$. Using further the weak lower semicontinuity of norms to pass to the limit in both the energy inequality \eqref{-271} and the uniform bounds \eqref{-66} we have proved the following proposition:

\begin{proposition}
\label{System on eta-level}
Let $\eta, \epsilon, \alpha > 0$, $\beta > \max\lbrace 4, \gamma \rbrace$ sufficiently large and let the assumptions of Theorem \ref{mainresult} be satisfied. Furthermore, let
\begin{align}
\rho _0 \in C^{2,\frac{1}{2}}\left(\overline{\Omega}\right), \quad \quad (\rho u)_0 \in C^2&\left(\overline{\Omega}\right), \quad \quad B_0 \in H^{2,2}(\Omega), \label{-342} \\
0 < \alpha \leq \rho_0 \leq \alpha^{-\frac{1}{2\beta}},\quad \quad \left. \nabla \rho_0 \cdot \operatorname{n} \right|_{\partial \Omega} = 0,& \quad \quad \operatorname{div} B_0 = 0,\quad \quad \left. B_0 \cdot \operatorname{n}\right|_{\partial \Omega} = 0. \label{-343}
\end{align}

Under these conditions there exist functions $X_\eta: \mathbb{R}^3 \rightarrow \mathbb{R}^3$, $u_\eta \in L^2(0,T;H_0^{1,2}(\Omega))$, $z_\eta, \tilde{z}_\eta \in L^\frac{4}{3} (Q)$ and 
\begin{align}
0 \leq \rho_\eta \in& \bigg\lbrace \psi \in L^\infty\left(0,T;L^\beta(\Omega)\right) \bigcap L^{r} \left(0,T; W^{1,r}(\Omega) \right) \bigcap L^{\tilde{r}}\left(0,T; W^{2,\tilde{r}}(\Omega) \right): \nonumber \\
&\partial_t \psi \in L^{\tilde{r}} \left(Q \right),\ \left. \nabla \psi \cdot \operatorname{n} \right|_{\partial \Omega} = 0 \bigg\rbrace, \label{-292} \\
B_\eta \in& \left\lbrace b \in L^\infty \left(0,T;L^2(\Omega)\right) \bigcap L^2\left(0,T;H^{2,2}_{\operatorname{div}}(\Omega)\right):\ \operatorname{curl}b = 0 \ \text{in } Q^s\left(S_\eta \right),\ b \cdot \operatorname{n}|_{\partial \Omega} = 0 \right\rbrace, \nonumber
\end{align}

for $r>2$, $\tilde{r}>1$ as in \eqref{-347} and $S_\eta= S_\eta (\cdot) = (X_\eta (\cdot; O))^\delta$, which satisfy
\begin{align}
\frac{dX_\eta(t;x)}{dt} = R_\delta \left[ u_\eta \right] \left( t, X_\eta(t;x) \right),\quad &\text{for } t \in [0,T], \quad \quad X_\eta(0;x) =  x,\quad \text{for } x \in \mathbb{R}^3, \nonumber \\
\partial_t \rho_\eta + \operatorname{div} \left( \rho_\eta u_\eta \right) &= \epsilon \Delta \rho_\eta \quad \text{a.e. in } Q \label{-296}
\end{align}

and
\begin{align}
-\int_0^T \int_\Omega \rho_\eta u_\eta \cdot \partial_t  \phi \ dxdt =& \int _0^T \int _\Omega \left( \rho_\eta u_\eta \otimes u_\eta \right): \mathbb{D} (\phi) + \left( a\rho_\eta^\gamma + \alpha \rho_\eta^\beta \right) \operatorname{div} \phi - 2\nu \left( \chi_\eta \right) \mathbb{D}(u_\eta):\mathbb{D} (\phi) \nonumber \\
&- \lambda \left( \chi_\eta \right) \operatorname{div} u_\eta \operatorname{div} \phi + \rho_\eta g\cdot \phi + \frac{1}{\mu}\left( \operatorname{curl}B_\eta \times B_\eta \right) \cdot \phi - \epsilon \tilde{z}_\eta \cdot \phi \nonumber \\
&- \epsilon \left( \nabla u_\eta \nabla \rho_\eta \right) \cdot \phi \ dxdt, \nonumber \\
- \int _0^T \int_\Omega B_\eta \cdot \partial_t b\ dxdt =& \int _0^T \int _\Omega \left[ -\frac{1}{\sigma \mu} \operatorname{curl} B_\eta + u_\eta \times B_\eta + \frac{1}{\sigma} J - \frac{\epsilon}{\mu ^2} z_\eta \right] \cdot \operatorname{curl} b \nonumber \\
&- \epsilon \operatorname{curl}\left( \operatorname{curl} B_\eta \right): \operatorname{curl}\left( \operatorname{curl} b\right)\ dxdt, \label{-298}
\end{align}

where $\chi_{\eta}(t,x) := \textbf{db}_{S_{\eta}(t)}(x)$, for any $\phi \in \mathcal{D}(Q)$ and any $b\in Y(S_\eta)$. Further, these functions satisfy the initial conditions
\begin{align}
\rho_\eta(0) = \rho_0,\quad \quad \left( \rho_\eta u_\eta \right)(0) = \left( \rho u \right)_0,\quad \quad B_\eta(0) = B_0, \nonumber
\end{align}

of which the latter identity can be understood in the sense of \eqref{-341}, as well as the energy inequality
\begin{align}
&\int_\Omega \frac{1}{2} \rho_\eta (t)\left|u_\eta (t)\right|^2 + \frac{a}{\gamma - 1} \rho_\eta ^\gamma(t) + \frac{\alpha}{\beta - 1}\rho_\eta ^\beta(t) + \frac{1}{2\mu} \left| B_\eta(t) \right|^2\ dx + \int_0^t \int_\Omega 2\nu \left(\chi_\eta \right) |\mathbb{D} \left(u_\eta \right)|^2 \nonumber \\
&+ \lambda \left(\chi_\eta \right)\left| \operatorname{div}u_\eta \right|^2 \nonumber + \alpha \epsilon \beta \rho_\eta^{\beta - 2}\left| \nabla \rho_\eta \right|^2 + \epsilon \left| \tilde{z}_\eta \right|^\frac{4}{3} + \frac{1}{\sigma \mu^2} \left| \operatorname{curl} B_\eta \right|^2 + \frac{\epsilon}{\mu} \left| \Delta B_\eta \right|^2 + \frac{\epsilon}{\mu^3} \left| z_\eta \right|^\frac{4}{3} \ dxd\tau \nonumber \\
\leq& \int_\Omega  \frac{1}{2} \rho_0\left|u_0\right|^2 + \frac{a}{\gamma - 1} \rho_0 ^\gamma + \frac{\alpha}{\beta - 1}\rho_0 ^\beta + \frac{1}{2\mu} \left| B_0 \right|^2\ dx +\int_{0}^t\int_\Omega \rho_\eta g \cdot u_\eta + \frac{1}{\sigma \mu} J \cdot \operatorname{curl} B_\eta \ dxd\tau \label{-300}
\end{align}

for almost all $t \in [0,T]$ and the estimate
\begin{align}
\epsilon \left\| \partial_t \rho_\eta \right\|_{L^{\tilde{r}}(Q)} + \epsilon^2 \left\| \Delta \rho_\eta \right\|_{L^{\tilde{r}}(Q)} + \epsilon \left\| \nabla \rho_\eta \right\|_{L^{r}(Q)} \leq c \label{-293}
\end{align}

with a constant $c>0$ independent of $\eta$ and $\epsilon$.

\end{proposition}

From this point on, the remainder of the proof of the main result is straight forward: In the mechanical part of the problem we can follow precisely the arguments from \cite[Sections 7--9]{feireisl}, the additional Lorentz force (c.f.\ \cite{sart}) and regularization term in the momentum equation do not cause any essential further difficulties. In the induction equation, each limit passage from now on can be carried out as in the incompressible case in \cite{incompressiblecase} and thus essentially as in Section \ref{inductionequationetalevel}. However, for the convenience of the reader, we will sketch the main arguments for the remaining three limit passages in the following sections.

{\centering \section{Limit passage in the penalization method} \par \label{etalimit}} 
We continue by passing to the limit with respect to $\eta \rightarrow 0$. Exactly as in the limit passage with respect to $n \rightarrow \infty$ in Section \ref{galerkinlimit} we can, due to the energy inequality \eqref{-300} and the uniform bound \eqref{-293}, extract suitable subsequences and find functions $z,\tilde{z} \in L^\frac{4}{3}(Q)$ and
\begin{align}
u \in& L^2(0,T;H_0^{1,2}(\Omega)), \label{-361} \\
0 \leq \rho \in& \bigg\lbrace \psi \in L^\infty\left(0,T;L^\beta(\Omega)\right) \bigcap L^{r} \left(0,T; W^{1,r}(\Omega) \right) \bigcap L^{\tilde{r}}\left(0,T; W^{2,\tilde{r}}(\Omega) \right): \nonumber \\
&\partial_t \psi \in L^{\tilde{r}} \left(Q \right),\ \left. \nabla \psi \cdot \operatorname{n} \right|_{\partial \Omega} = 0 \bigg\rbrace \subset C \left([0,T];L^2(\Omega) \right), \label{-362} \\
B \in& \left\lbrace b \in L^\infty \left(0,T;L^2(\Omega)\right) \bigcap L^2\left(0,T;H^{2,2}_{\operatorname{div}}(\Omega)\right):\ \operatorname{curl}b = 0 \ \text{in } Q^s(S),\ b \cdot \text{n}|_{\partial \Omega} = 0 \right\rbrace, \label{-290}
\end{align}

such that
\begin{align}
\rho_\eta \rightarrow \rho \ \ \ &\text{in } L^\beta \left(Q \right), \quad \quad \quad &\rho_\eta& \rightarrow \rho \quad \ \ \ \text{in } L^2 \left(0,T;H^{1,2}(\Omega)\right), \label{-100} \\
\rho_\eta \rightharpoonup \rho \ \ \ &\text{in } L^{\tilde{r}}\left(0,T;W^{2,\tilde{r}}(\Omega) \right), \quad \quad \quad &\partial_t \rho_\eta & \rightharpoonup \partial_t \rho \ \ \ \ \text{in } L^{\tilde{r}}\left(Q \right), \label{-99} \\
u_\eta \rightharpoonup u \ \ \ &\text{in } L^2 \left(0,T;H^{1,2}(\Omega)\right), \quad \quad \quad &B_\eta& \buildrel\ast\over\rightharpoonup B \ \ \quad \text{in } L^\infty \left(0,T;L^2(\Omega)\right), \label{-97} \\
B_\eta \rightharpoonup B \ \ \ &\text{in } L^2 \left(0,T;H^{2,2}(\Omega)\right), \quad \quad \quad &z_\eta& \rightharpoonup z \quad \quad \text{in } L^\frac{4}{3}\left( Q \right), \nonumber \\
\tilde{z}_\eta \rightharpoonup \tilde{z} \ \quad &\text{in } L^\frac{4}{3}\left( Q \right). \nonumber
\end{align}

The boundary conditions of the limit functions in \eqref{-361}--\eqref{-290} follow directly from the corresponding boundary conditions on the $\eta$-level, see Proposition \ref{System on eta-level}. The set-valued function $S$ in \eqref{-290} is defined by $S:=S(\cdot) := (X(\cdot;O))^\delta$ where $X$, given by
\begin{align}
X_\eta \rightarrow X \quad \text{in } C\left([0,T];C_{\text{loc}}\left( \mathbb{R}^3 \right) \right), \nonumber
\end{align}

denotes the solution to the initial value problem
\begin{align}
\frac{dX(t;x) }{dt} = R_\delta [u] \left( t, X (t;x) \right),\quad X (0;x)=x \quad \quad \forall x \in \mathbb{R}^3, \label{-349}
\end{align}

c.f.\ Lemma \ref{convergence}. In particular, this lemma also implies that
\begin{align}
\chi_\eta \rightarrow \chi := \textbf{db}_{S(\cdot)}(\cdot) \quad \text{in } C\left([0,T];C_{\text{loc}}\left( \mathbb{R}^3 \right) \right). \nonumber
\end{align}

\subsection{Continuity equation}

Making use of the convergences \eqref{-100}--\eqref{-97} of $\rho_\eta$ and $u_\eta$, we can pass to the limit in the continuity equation \eqref{-296} and ensure that
\begin{align}
\partial_t \rho + \operatorname{div} \left( \rho u \right) =& \epsilon \Delta \rho \quad \text{a.e. in } Q,\quad \rho(0) = \rho_0. \label{-103}
\end{align}

Moreover, this pointwise identity can be renormalized by multiplying it by $\zeta'(\rho_\epsilon)$ for an arbitrary convex function $\zeta \in C^2([0,+\infty))$. Since $\zeta '' \geq 0$, this yields
\begin{align}
\partial_t \zeta \left(\rho \right) + \operatorname{div} \left( \zeta \left(\rho \right)u \right) + \left[ \zeta' \left(\rho \right)\rho - \zeta \left(\rho \right) \right] \operatorname{div} u - \epsilon \Delta \zeta \left(\rho \right) = - \epsilon \zeta'' \left(\rho \right)\left| \nabla \rho \right|^2 \leq 0 \label{-152}
\end{align}

almost everywhere in $Q$. This relation will turn out useful in the limit passage with respect to $\epsilon \rightarrow 0$ in Section \ref{epsilonlimit}.

\subsection{Induction equation}

For the limit passage in the induction equation we can argue exactly as in the limit passage with respect to $n \rightarrow \infty$ in Section \ref{inductionequationetalevel} to show strong convergence of $B_\eta$ in the fluid domain. Hence, we can pass to the limit in \eqref{-298} and obtain the identity
\begin{align}
-\int_0^T \int_\Omega B \cdot \partial_t b\ dxdt =& \int_0^T\int_\Omega \left[ -\frac{1}{\sigma \mu} \operatorname{curl} B + u \times B - \frac{\epsilon}{\mu ^2} z + \frac{1}{\sigma}J \right] \cdot \operatorname{curl} b \nonumber \\
&- \epsilon \operatorname{curl}\left(\operatorname{curl} B\right) \cdot \operatorname{curl}\left(\operatorname{curl} b\right)\ dxdt \quad \quad \forall b \in Y(S). \label{-109}
\end{align}

Moreover, the initial condition $B(0)=B_0$ also follows as in Section \ref{inductionequationetalevel}.

\subsection{Momentum equation and compatibility of the velocity field}\label{highviscositylimitmomentum}

From the uniform bounds given by the energy inequality \eqref{-300} we further infer the existence of $z_5 \in L^\frac{6}{5}(Q)$ such that, for a chosen subsequence,
\begin{align}
\rho_\eta u_\eta \otimes \eta \rightharpoonup z_5 \quad \text{in } L^\frac{6}{5}(Q). \nonumber
\end{align}

For the limit passage in the momentum equation we need to identify $z_5$: The choice of test functions $\phi \in T(S)$, which satisfy $\mathbb{D}(\phi) = 0$ in a neighbourhood of $\overline{Q}^s(S)$, allows us to control the variable viscosity coefficients $\nu (\chi_\eta)$ and $\lambda (\chi_\eta)$ in the momentum equation \eqref{-298}, since these remain bounded in the fluid region according to their definition in \eqref{-322}, \eqref{-52}. This enables us to deduce strong convergence of the momentum function $\rho_\eta u_\eta$ in the fluid domain similarly to the strong convergence \eqref{-282} of the magnetic induction in the Galerkin limit. Indeed, we fix an arbitrary interval $I \subset (0,T)$ and an arbitrary ball $U \subset \Omega$ such that $\overline{I \times U} \subset Q^f(S)$ and deduce from the momentum equation, for any $\Phi \in \mathcal{D}(U)$, the dual estimate
\begin{align}
\left\| \partial_t \int_U \rho_\eta u_\eta \cdot \Phi dx \right\|_{L^{\min\left(\frac{6}{5}, \frac{5\beta - 3}{4\beta}\right)}(I)} \leq c \nonumber
\end{align}

for a constant $c>0$ depending on $\Phi$ but not on $\eta$. From the Arzelá-Ascoli theorem it follows that
\begin{align}
\rho_\eta u_\eta \rightarrow \rho u \quad \text{in } C_\text{weak}\left( \overline{I};L^\frac{2\beta}{\beta + 1}(U) \right) \quad \text{and hence in } L^2 \left( I;H^{-1,2}(U) \right), \nonumber
\end{align}

which implies that
\begin{align}
z_5 : \mathbb{D}(\phi) = (\rho u \otimes u ): \mathbb{D}(\phi) \quad \text{a.e. in } Q. \nonumber
\end{align}

for all test functions $\phi \in T(S)$. Letting $\eta$ tend to $0$ in \eqref{-298} we thus obtain
\begin{align}
-\int_{0}^{T} \int_\Omega \rho u \cdot \partial_t \phi\ dxdt =& \int _{0}^{T} \int _\Omega \left( \rho u \otimes u \right): \mathbb{D} (\phi) + \left( a\rho^\gamma + \alpha \rho^\beta \right) \operatorname{div} \phi - 2\nu \mathbb{D}\left(u \right):\mathbb{D} (\phi) - \lambda \operatorname{div} u\operatorname{div} \phi \nonumber \\
&+ \rho g\cdot \phi + \frac{1}{\mu}\left( \text{curl}B \times B \right) \cdot \phi - \epsilon \tilde{z} \cdot \phi - \epsilon \left( \nabla u \nabla \rho \right) \cdot \phi\ dxdt \quad \quad \forall \phi \in T(S). \label{-112}
\end{align}

Moreover, since $\nu (\chi_\eta)$ and $\lambda (\chi_\eta)$ blow up in the solid part of the domain, the energy inequality \eqref{-300} shows that
\begin{align}
\mathbb{D}(u)=0 \quad \text{a.e. in } Q^s(S). \nonumber
\end{align}

Hence, there are rigid velocity fields $u^{s^i}$ which coincide with $u$ almost everywhere in the $\delta$-neighbourhoods $S^i(t):= (X(t;O^i))^\delta$ of the sets $X(t;O^i)$. Consequently, due to the property \eqref{-348} of the regularized velocity field $R_\delta [u]$, we can replace $R_\delta[u]$ in the initial value problem \eqref{-349} by $u^{s^i}$ for $x \in O^i$. The combination of the latter two conditions at first yields compatibility (c.f.\ \eqref{-344}, \eqref{-345}) of $u$ with the system $\lbrace O^i, X^i \rbrace_{i=1}^m$, where each $X^i(t)$ denotes an isometry which coincides with $X(t)$ in $O^i$. However, the fact that each $X^i(t)$ is an isometry implies that
\begin{align}
S^i(t) = \left(X^i\left(t;O^i \right)\right)^\delta = X^i\left(t; S^i_0 \right) \quad \quad \text{and thus } \quad \quad u(t,\cdot) = u^{s^i}(t,\cdot) \quad \text{a.e. in } X^i\left(t; S^i \right), \nonumber
\end{align}

Consequently, we infer that $u$ is even compatible with $\lbrace S_0^i, X^i \rbrace_{i=1}^m$. Finally, the initial condition $(\rho u)(0) = (\rho u)_0$, in the sense of \eqref{-341}, follows by similar arguments as the initial condition for the magnetic induction in Section \ref{inductionequationetalevel}.

\subsection{Energy inequality}\label{energyinequalityetalimit}

We drop, among other non-negative terms, the variable parts of the viscosity coefficients on the left-hand side of the energy inequality \eqref{-300} and pass to the limit to see that
\begin{align}
&\int_\Omega \frac{1}{2} \rho (t)\left|u (t)\right|^2 + \frac{a}{\gamma - 1} \rho ^\gamma(t) + \frac{\alpha}{\beta - 1}\rho ^\beta(t) + \frac{1}{2\mu} \left| B(t) \right|^2\ dx + \int_0^t \int_\Omega 2\nu |\mathbb{D} (u)|^2 + \lambda \left| \operatorname{div}u \right|^2 + \epsilon \left| \tilde{z} \right|^\frac{4}{3} \nonumber \\
&+ \frac{1}{\sigma \mu^2} \left| \operatorname{curl} B \right|^2 + \frac{\epsilon}{\mu} \left| \Delta B \right|^2 + \frac{\epsilon}{\mu^3} \left| z \right|^\frac{4}{3}\ dxd\tau \nonumber \\
\leq& \int_\Omega  \frac{1}{2} \rho_0\left|u_0\right|^2 + \frac{a}{\gamma - 1} \rho_0 ^\gamma + \frac{\alpha}{\beta - 1}\rho_0 ^\beta + \frac{1}{2\mu} \left| B_0 \right|^2\ dx +\int_{0}^t\int_\Omega \rho g \cdot u + \frac{1}{\sigma \mu} J \cdot \operatorname{curl} B \ dxd\tau \label{-116}
\end{align}

for almost all $t \in [0,T]$.

{\centering \section{Limit passage in the regularization terms} \label{epsilonlimit} \par }

The next step is the limit passage with respect to $\epsilon \rightarrow 0$. The energy inequality \eqref{-116} yields the bounds needed for Corollary \ref{convergenceandcompatibility}, which implies the existence of isometries $X^i(t):\mathbb{R}^3 \rightarrow \mathbb{R}^3$ such that
\begin{align}
X^i_\epsilon \rightarrow X^i \quad &\text{in } C\left([0,T];C_{\text{loc}}\left( \mathbb{R}^3 \right) \right). \nonumber
\end{align}

We write $X:[0,T]\times S_0 \rightarrow \mathbb{R}^3$,\ $X(t)|_{S_0^i}:= X^i(t)$ and $S=S(\cdot) := X(\cdot; S_0)$. Further, testing the continuity equation \eqref{-103} by $\rho _\epsilon$, we see that
\begin{align}
\epsilon^\frac{1}{2} \left\| \nabla \rho_\epsilon \right\|_{L^2(Q)} \leq c \nonumber
\end{align}

for a constant $c>0$ independent of $\epsilon$. This, together with the energy inequality \eqref{-116}, yields the existence of functions $0 \leq \rho \in L^\infty(0,T;L^\beta (\Omega))$ and
\begin{align}
u \in& \left\lbrace \phi \in L^2\left(0,T;H_0^{1,2}\left(\Omega \right) \right):\ \mathbb{D}(\phi) = 0\ \text{in } Q^s(S) \right\rbrace, \label{-363} \\
B \in& \left\lbrace b \in L^\infty \left(0,T;L^2(\Omega)\right) \bigcap L^2\left(0,T;H^{1,2}_{\operatorname{div}}(\Omega)\right):\ \operatorname{curl}b = 0 \ \text{in } Q^s(S),\ b \cdot \text{n}|_{\partial \Omega} = 0 \right\rbrace \label{-364}
\end{align}

such that for certain extracted subsequences it holds
\begin{align}
\rho_\epsilon \buildrel\ast\over\rightharpoonup \rho \ \ \ &\text{in } L^\infty \left(0,T;L^\beta(\Omega)\right), \quad \quad \quad &u_\epsilon& \rightharpoonup u \ \ \ \ \text{in } L^2 \left(0,T;H^{1,2}(\Omega)\right), \nonumber \\
B_\epsilon \buildrel\ast\over\rightharpoonup B \ \ \ &\text{in } L^\infty \left(0,T;L^2(\Omega)\right), \ \ \ \ \ \ &B_\epsilon& \rightharpoonup B \ \ \ \text{in } L^2 \left(0,T;H^{1,2}(\Omega)\right), \nonumber \\
\epsilon \nabla \rho _\epsilon,\ \epsilon \Delta B_\epsilon \rightarrow 0 \quad &\text{in } L^2\left(Q\right), \quad \quad &\epsilon z_\epsilon,\ \epsilon \tilde{z}_\epsilon& \rightarrow 0 \quad \ \text{in } L^\frac{4}{3} \left(Q\right). \label{-135}
\end{align}

The boundary conditions of the limit functions in \eqref{-363} and \eqref{-364} follow directly from the corresponding boundary conditions of the velocity field and the magnetic induction in \eqref{-361} and $\eqref{-290}$ on the $\epsilon$-level. Moreover, the velocity field $u$ is compatible with the system $\lbrace S^i_0, X^i \rbrace_{i=1}^m$, c.f.\ Corollary \ref{convergenceandcompatibility}.\\

\subsection{Continuity equation} \label{contequationepsilonlimit}

Similarly to the strong convergence \eqref{-282} of the magnetic induction, we deduce from the continuity equation \eqref{-103} that $\rho_\epsilon$ even converges to $\rho$ in $C_ \text{weak}([0,T];L^\beta (\Omega))$. This, together with the vanishing artificial viscosity term, c.f.\ \eqref{-135}, is sufficient to pass to the limit in \eqref{-103} and to obtain
\begin{align}
\partial_t \rho + \operatorname{div} \left( \rho u \right) = 0 \quad \text{in } \mathcal{D}'\left(Q \right). \label{-142}
\end{align}

In fact, as at this stage of the approximation it holds $\rho \in L^2(Q)$, $u \in L^2(0,T;H_0^{1,2}(\Omega))$, we can use the regularization procedure by DiPerna and Lions, c.f.\ \cite[Lemma 6.8, Lemma 6.9]{novotnystraskraba}, to see that $\rho$ and $u$, extended by $0$ outside of $\Omega$, even satisfy the renormalized continuity equation \eqref{-215}, \eqref{-216}. This in turn implies that $\rho \in C([0,T];L^1(\Omega))$, c.f.\ \cite[Lemma 6.15]{novotnystraskraba}, and $\rho$ satisfies the initial condition $\rho(0)=\rho_0$.

\subsection{Induction equation}

The regularization terms in the induction equation vanish as $\epsilon$ tends to $0$ according to the convergences \eqref{-135}. Apart from that we can argue exactly as in the Galerkin limit in Section \ref{inductionequationetalevel} to pass to the limit with respect to $\epsilon \rightarrow 0$ in \eqref{-109} and to infer that
\begin{align}
-\int_0^T \int_\Omega B \cdot \partial_t b\ dxdt = \int_0^T\int_\Omega \left[ -\frac{1}{\sigma \mu} \operatorname{curl} B + u \times B + \frac{1}{\sigma}J \right] \cdot \operatorname{curl} b\ dxdt \quad \quad \forall b \in Y(S). \label{-136}
\end{align}

\subsection{Momentum equation} \label{momentumequationepsilonlimit}

In order to pass to the limit in the pressure terms, we first consider an arbitrary compact set $K \subset
Q^f(S)$. Denoting by $\mathcal{B}_\Omega$ the Bogovskii operator in $\Omega$ (c.f.\ \cite[Section 3.3.1.2]{novotnystraskraba}), we test the momentum equation \eqref{-112} by
\begin{align}
\phi_\epsilon (t,x) := \Phi(t,x) \mathcal{B}_\Omega \left[ \rho_\epsilon (t,\cdot) - \frac{1}{|\Omega|} \int_\Omega \rho_\epsilon (t,y)\ dy \right](t,x), \label{-350}
\end{align}

where $\Phi \in \mathcal{D}(Q^f(S))$ is a cut-off function equal to $1$ in $K$. This procedure leads to a bound of $\rho_\epsilon$ in $L^{\beta + 1}(K)$ uniformly in $\epsilon$, c.f.\ \cite[Lemma 8.1]{feireisl} and the references therein. These bounds in turn allow us to find $z_6 \in L^\frac{\gamma + 1}{\gamma}(K)$, $z_7 \in L^\frac{\beta + 1}{\beta}(K)$ such that
\begin{align}
\rho_\epsilon^\gamma \rightharpoonup z_6 \quad \text{in } L^\frac{\gamma + 1}{\gamma} \left(K\right), \quad \quad \rho_\epsilon^\beta \rightharpoonup z_7 \quad \text{in } L^\frac{\beta + 1}{\beta} \left(K \right). \nonumber
\end{align}

With the aim of identifying these limit functions we set, for arbitrary $\tilde{\Phi} \in \mathcal{D}(Q^f(S))$,
\begin{align}
\tilde{\phi}_\epsilon (t,x) := \tilde{\Phi}(t,x) \left( \nabla \Delta^{-1} \right)\left[ \rho_\epsilon (t,\cdot) \right](t,x),\quad \quad \tilde{\phi} (t,x) := \tilde{\Phi}(t,x) \left( \nabla \Delta^{-1} \right)\left[ \rho(t,\cdot) \right](t,x), \label{-351}
\end{align}

where $\Delta^{-1}$ denotes the inverse Laplacian on $\mathbb{R}^3$, c.f.\ \cite[Section 10.16]{singularlimits}. We compare the momentum equation \eqref{-112} on the $\epsilon$-level, tested by $\tilde{\phi}_\epsilon$, to a corresponding limit identity, tested by $\tilde{\phi}$. This enables us to deduce the effective viscous flux identity
\begin{align}
(\lambda + 2\nu) \lim_{\epsilon \rightarrow 0} \int_{0}^T \int_{\Omega} \tilde{\Phi} \left( \rho_\epsilon \operatorname{div} u_\epsilon - \rho \operatorname{div} u \right) \ dxdt = \lim_{\epsilon \rightarrow 0} \int_{0}^T \int_{\Omega} \tilde{\Phi} \left( \left[ a\rho_\epsilon^\gamma + \alpha \rho_\epsilon^\beta \right] \rho_\epsilon - \left[ a\rho_\epsilon^{\gamma } + \alpha \rho_\epsilon^{\beta } \right] \rho \right)\ dxdt \nonumber
\end{align}

for all $\tilde{\Phi} \in \mathcal{D}(Q^f(S))$, c.f.\ \cite[Lemma 8.2]{feireisl} and the references therein. Moreover, after a density argument, we can consider the choice $\zeta(s) = s \ln(s)$ in both the renormalized continuity equations \eqref{-152} on the $\epsilon$-level and \eqref{-215} in the limit. A comparison between the resulting identities then leads us to
\begin{align}
 \int_\Omega \zeta \left( \rho (\tau) \right)\ dx - \lim_{\epsilon \rightarrow 0} \int_\Omega \zeta \left( \rho_\epsilon (\tau) \right)\ dx \geq \lim_{\epsilon \rightarrow 0} \int_0^\tau \int_\Omega \rho_\epsilon \operatorname{div} u_\epsilon \ dxdt - \int_0^\tau \int_\Omega \rho \operatorname{div} u \ dxdt \geq 0 \label{-157}
\end{align}

for $\tau \in [0,T]$, where the last inequality follows from the effective viscous flux identity and the monotonicity of the mapping $s \mapsto as^\gamma + \alpha s^\beta$, as well as the fact that $\operatorname{div}u =0$ in $Q^s(S)$. Due to the strict convexity of $\zeta$ this estimate implies pointwise convergence of $\rho_\epsilon$ in $(0,T)\times \Omega$, c.f.\ \cite[Theorem 10.20]{singularlimits}, and hence $z_6 = \rho^\gamma$, $z_7 = \rho ^\beta$ almost everywhere in $Q^f(S)$. In the remaining terms of the momentum equation \eqref{-112} we can pass to the limit as during the past limit passages. We end up with
\begin{align}
- \int_0^T \int_\Omega \rho u \cdot \partial_t \phi \ dxdt =& \int _0^T \int _\Omega \left( \rho u \otimes u \right): \mathbb{D} (\phi) + \left( a\rho^\gamma + \alpha \rho ^\beta \right) \operatorname{div} \phi - 2\nu \mathbb{D}(u):\mathbb{D} (\phi) \nonumber \\
&- \lambda \operatorname{div} u \operatorname{div} \phi + \rho g\cdot \phi + \frac{1}{\mu}\left( \operatorname{curl}B \times B \right) \cdot \phi \ dxdt \quad \quad \forall \phi \in T(S). \label{-312}
\end{align}

\subsection{Energy inequality}

Neglecting the regularization terms on the left-hand side of the energy inequality \eqref{-116} on the $\epsilon$-level, we can pass to the limit with respect to $\epsilon \rightarrow 0$ and obtain
\begin{align}
&\int_\Omega \frac{1}{2} \rho (t)\left|u (t)\right|^2 + \frac{a}{\gamma - 1} \rho ^\gamma(t) + \frac{\alpha}{\beta - 1}\rho ^\beta(t) + \frac{1}{2\mu} \left| B(t) \right|^2\ dx + \int_0^t \int_\Omega 2\nu |\mathbb{D} (u)|^2 + \lambda \left| \operatorname{div}u \right|^2 \nonumber \\
&+ \frac{1}{\sigma \mu^2} \left| \operatorname{curl} B \right|^2\ dxd\tau \nonumber \\
\leq& \int_\Omega  \frac{1}{2} \rho_0\left|u_0\right|^2 + \frac{a}{\gamma - 1} \rho_0 ^\gamma + \frac{\alpha}{\beta - 1}\rho_0 ^\beta + \frac{1}{2\mu} \left| B_0 \right|^2\ dx +\int_{0}^t\int_\Omega \rho g \cdot u + \frac{1}{\sigma \mu} J \cdot \operatorname{curl} B \ dxd\tau \label{-167}
\end{align}

for almost every $t \in [0,T]$.

{\centering \section{Limit passage in the artificial pressure} \label{alphalimit} \par }

Finally it remains to pass to the limit with respect to $\alpha \rightarrow 0$. We now consider initial data $\rho_0$, $(\rho u)_0$ and $B_0$ as in Theorem \ref{mainresult} and construct - c.f.\ \cite[Section 4]{feireisl2} - the initial data $\rho_{0,\alpha}$, $(\rho u)_{0,\alpha}$ and $B_{0,\alpha}$ on the $\alpha$-level (c.f.\ \eqref{-342}, \eqref{-343}) in such a way that
\begin{align}
\rho_{0,\alpha} \rightarrow \rho_0 \quad \ \ \ &\text{in } L^\gamma (\Omega),\quad \quad &\alpha \rho_{0,\alpha}^\beta& \rightarrow 0 \quad \quad \quad \quad \text{in } L^1(\Omega), \nonumber \\
\left(\rho u \right)_{0,\alpha} \rightarrow \left( \rho u \right)_0 \quad &\text{in } L^1 (\Omega),\quad \quad &\frac{\left|\left(\rho u \right)_{0,\alpha}\right|^2}{\rho_{0,\alpha}}& \rightarrow \frac{\left|\left(\rho u \right)_{0}\right|^2}{\rho_{0}} \quad \text{in } L^1(\Omega), \nonumber \\
B_{0,\alpha} \rightarrow B_0 \ \ \ \quad \ &\text{in } L^2(\Omega). && \nonumber
\end{align}

Since the energy inequality \eqref{-116} provides the conditions for Corollary \ref{convergenceandcompatibility}, we obtain the existence of isometries $X^i(t):\mathbb{R}^3 \rightarrow \mathbb{R}^3$ such that
\begin{align}
X^i_\alpha \rightarrow X^i \quad &\text{in } C\left([0,T];C_{\text{loc}}\left( \mathbb{R}^3 \right) \right). \nonumber
\end{align}

We set 
\begin{align}
X:[0,T]\times S_0 \rightarrow \mathbb{R}^3,\quad X(t)|_{S_0^i}:= X^i(t) \label{-327}
\end{align}

and $S=S(\cdot) := X(\cdot; S_0)$. Then from the energy inequality \eqref{-167} we obtain the existence of functions
\begin{align}
0 \leq \rho \in& L^\infty(0,T;L^\gamma (\Omega)), \label{-313} \\
u \in& \left\lbrace \phi \in L^2\left(0,T;H_0^{1,2}\left(\Omega \right) \right):\ \mathbb{D}(\phi) = 0\ \text{in } Q^s(S) \right\rbrace \label{-309} \\
B \in& \left\lbrace b \in L^\infty \left(0,T;L^2(\Omega)\right) \bigcap L^2\left(0,T;H^{1,2}_{\operatorname{div}}(\Omega)\right):\ \operatorname{curl}b = 0 \ \text{in } Q^s(S),\ b \cdot \text{n}|_{\partial \Omega} = 0 \right\rbrace \label{-310}
\end{align}

such that for suitable subsequences it holds
\begin{align}
\rho_\alpha \buildrel\ast\over\rightharpoonup \rho \ \ \ &\text{in } L^\infty \left(0,T;L^\gamma(\Omega)\right), \quad \quad \quad &u_\alpha& \rightharpoonup u \ \ \ \ \text{in } L^2 \left(0,T;H^{1,2}(\Omega)\right), \nonumber \\
B_\alpha \buildrel\ast\over\rightharpoonup B \ \ \ &\text{in } L^\infty \left(0,T;L^2(\Omega)\right), \ \ \ \ \ \ &B_\alpha& \rightharpoonup B \ \ \ \text{in } L^2 \left(0,T;H^{1,2}(\Omega)\right). \nonumber
\end{align}

The boundary conditions of the limit functions in \eqref{-309} and \eqref{-310} follow directly from the corresponding boundary conditions of the velocity field and the magnetic induction in \eqref{-363} and $\eqref{-364}$ on the $\alpha$-level. Moreover, again due to Corollary \ref{convergenceandcompatibility},
\begin{align}
u\quad \text{is compatible with the family}\quad \lbrace S^i_0, X^i \rbrace_{i=1}^m. \label{-314}
\end{align}

\subsection{Continuity equation}

After using the continuity equation \eqref{-142} to deduce convergence of $\rho_\alpha$ in $C_\text{weak}([0,T];L^\gamma (\Omega))$, we pass to the limit in \eqref{-142} and obtain
\begin{align}
\partial_t \rho + \operatorname{div} \left( \rho u \right) = 0 \quad \text{in } \mathcal{D}'\left((0,T) \times \mathbb{R}^3 \right). \label{-183}
\end{align}

The proof of the renormalized continuity equation however needs to be postponed to Section \ref{momequationalphalimit} below, since at this stage $\rho$ does not have the $L^2(Q)$-regularity required for the regularization technique by DiPerna and Lions anymore.

\subsection{Induction equation}
For the limit passage in the induction equation \eqref{-136} we argue exactly as in Section \ref{inductionequationetalevel} and end up with
\begin{align}
-\int_0^T \int_\Omega B \cdot \partial_t b\ dxdt = \int_0^T\int_\Omega \left[ -\frac{1}{\sigma \mu} \operatorname{curl} B + u \times B + \frac{1}{\sigma}J \right] \cdot \operatorname{curl} b\ dxdt \quad \quad \forall \phi \in Y(S). \label{-186}
\end{align}

\subsection{Momentum equation} \label{momequationalphalimit}

For the limit passage in the pressure terms the strategy used during the limit passage with respect to $\epsilon \rightarrow 0$ in Section \ref{momentumequationepsilonlimit} needs to be modified to make up for the lower integrability of the density; the main ideas however remain the same. First we test the momentum equation \eqref{-312} by functions of the form \eqref{-350} with the density replaced by (a cut-off and smoothened version of) $\rho_\alpha^\theta$, $\theta > 0$. Choosing $\theta > 0$ sufficiently small, we find that, for any compact $K \subset Q^f(S)$, $\rho_\alpha$ and $\alpha^\frac{1}{\beta + \theta}\rho_\alpha$ are bounded uniformly in $L^{\gamma + \theta}(K)$ and $L^{\beta + \theta}(K)$, respectively, c.f.\ \cite[Section 4.1]{feireisl2}, \cite[Proposition 2.3]{feireisl3}. In particular, there exists $z_8 \in L^\frac{\gamma + \theta}{\gamma}(K)$ such that, after the extraction of a subsequence,
\begin{align}
\rho_\alpha^\gamma \rightharpoonup z_8\quad \text{in } L^\frac{\gamma + \theta}{\gamma} \left(K\right), \quad \quad \alpha \rho_\alpha^\beta \rightharpoonup 0 \quad \text{in } L^\frac{\beta + \theta}{\beta} \left(K \right). \nonumber
\end{align}

In order to identify $z_7$, we again need to show strong convergence of $\rho_\alpha$. To this end we use the notion of the oscillation defect measure
\begin{align}
\textbf{osc}_{\gamma + 1}\left[\rho_\alpha \rightarrow \rho \right]\left(O\right) := \sup_{k \geq 1} \left[ \limsup_{\alpha \rightarrow 0} \int_O \left| T_k \left(\rho_\alpha \right) - T_k\left( \rho \right) \right|^{\gamma + 1}\ dxdt \right] \nonumber
\end{align}

for measurable sets $O \subset (0,T)\times \mathbb{R}^3$ and a concave cut-off function $T_k \in C^\infty([0,\infty))$, $k \in \mathbb{N}$, coinciding with the identity function on $[0,k]$ and with $2k$ on $[3k,\infty)$. The proof of the pointwise convergence of $\rho_\alpha$ can be divided into three main steps, each of which consists of showing one of the following three relations, respectively:
\begin{itemize}
\item[(i)] The effective viscous flux identity
\begin{align}
(\lambda + 2\nu)& \lim_{\alpha \rightarrow 0} \int_0^T \int_\Omega \phi \left( T_k \left( \rho_\alpha \right) \operatorname{div}u_\alpha - \underline{T_k(\rho)} \operatorname{div}u_\alpha \right)\ dxdt \nonumber \\
=& \lim_{\alpha \rightarrow 0} \int_0^T \int_\Omega \phi \left( a \rho_\alpha^\gamma T_k \left( \rho_\alpha \right) - a \rho_\alpha^\gamma \underline{T_k(\rho)} \right)\ dxdt, \label{-197}
\end{align}

where $\underline{T_k(\rho)}$ denotes a weak $L^1$-limit of $T_k(\rho_\alpha)$, holds true for any $\phi \in \mathcal{D}(Q^f(S))$.

\item[(ii)] The oscillation defect measure is bounded on $(0,T)\times \mathbb{R}^3$,
\begin{align}
\textbf{osc}_{\gamma + 1}\left[\rho_\alpha \rightarrow \rho \right]\left((0,T) \times \mathbb{R}^3 \right) < +\infty. \label{-198}
\end{align}
\item[(iii)] The renormalized continuity equation \eqref{-215}, \eqref{-216} is satisfied by $\rho$ and $u$.
\end{itemize}

The effective viscous flux identity \eqref{-197} can be shown by a comparison between the momentum equation \eqref{-312} on the $\alpha$-level and a corresponding limit identity, tested by suitably modified variants of the functions $\tilde{\phi}_\epsilon$ and $\tilde{\phi}$ in \eqref{-351} with the density replaced by $T_k(\rho_\alpha)$ and $\underline{T_k(\rho)}$, respectively. The details, in the case without rigid bodies, are given e.g.\ in \cite[Section 4.3]{feireisl2}, the adjustment to the fluid-structure case poses no further difficulties. The proof of the bound \eqref{-198} of the oscillation defect measure is split up into an estimate on $Q^s(S)$ and an estimate on $Q^f(S)$. From the representation of the density in the solid region in Lemma \ref{compatibility} (ii) it follows that $\rho_\alpha \rightarrow \rho$ in $L^1(K)$ for compact sets $K \subset Q^s(S)$ and thus $\textbf{osc}_{\gamma + 1}\left[\rho_\alpha \rightarrow \rho \right](Q^s(S))=0$. In the fluid region the bound is achieved, under exploitation of the effective viscous flux identity \eqref{-197}, by the same arguments as in the all-fluid case, c.f.\ \cite[Proposition 6.1]{feireislcompactness}. Finally, the renormalized continuity equation in the limit is also obtained exactly as in the all-fluid case, c.f.\ \cite[Proposition 7.1]{feireislcompactness}: the idea is to pass to the limit in the renormalized continuity equation \eqref{-215} on the $\alpha$-level for the choice $\zeta = T_k$. Thanks to the boundedness of $T_k$, the regularization technique by DiPerna and Lions (c.f.\ \cite[Lemma 6.9]{novotnystraskraba}) can be applied to the limit identity. Letting $k \rightarrow \infty$ and exploiting the bound \eqref{-198} of the oscillation defect measure, we then obtain the renormalized continuity equation \eqref{-215}, \eqref{-216} also for $\rho$ and $u$. \\

Having shown the relations (i)--(iii) we now obtain strong convergence of $\rho_\alpha$. Indeed, similarly as in the corresponding relation \eqref{-157} in the $\epsilon$-limit and under exploitation of the concavity of $T_k$, we see that the left-hand side of the effective viscous flux identity \eqref{-197} is non-negative. This, in combination with the bound \eqref{-198} of the oscillation defect measure and a comparison between the renormalized continuity equations on the $\alpha$-level and in the limit, yields, similarly to the first inequality in \eqref{-157},
\begin{align}
 \int_\Omega \zeta \left( \rho (\tau) \right)\ dx - \lim_{\alpha \rightarrow 0} \int_\Omega \zeta \left( \rho_\alpha (\tau) \right)\ dx \geq 0 \nonumber
\end{align}

for $\tau \in [0,T]$ and $\zeta (s) := s \ln (s)$. As in Section \ref{momentumequationepsilonlimit}, this inequality implies pointwise convergence of $\rho_\alpha$ in $(0,T)\times \Omega$ and therefore $z_8 = \rho^\gamma$ almost everywhere in $Q^f(S)$. In the remaining terms of the momentum equations \eqref{-312} we may pass to the limit as in the previous limit passages and obtain
\begin{align}
- \int_0^T \int_\Omega \rho u \cdot \partial_t \phi \ dxdt =& \int _0^T \int _\Omega \left( \rho u \otimes u \right): \mathbb{D} (\phi) + a\rho^\gamma \operatorname{div} \phi - 2\nu \mathbb{D}(u):\mathbb{D} (\phi) \nonumber \\
&- \lambda \operatorname{div} u \operatorname{div} \phi + \rho g\cdot \phi + \frac{1}{\mu}\left( \operatorname{curl}B \times B \right) \cdot \phi \ dxdt \quad \forall \phi \in T(S). \label{-224}
\end{align}

We are now in the position to conclude the proof of the main result.

\subsection{Proof of the main result}

The function $X$ in \eqref{-159} is defined in \eqref{-327}. The properties of $\rho$, $u$ and $B$ in \eqref{-211}--\eqref{-212}, except for the continuity of $\rho$ in time, are shown in \eqref{-313}--\eqref{-310}. The continuity equation \eqref{-214} and its renormalization \eqref{-215}, \eqref{-216} are derived in \eqref{-183} and the relation (iii) in Section \ref{momequationalphalimit}, respectively. In particular, $\rho$ as a renormalized solution to the continuity equation satisfies $\rho \in C([0,T];L^1(\Omega))$, c.f.\ \cite[Proposition 4.3]{feireisldynamics}, which concludes the proof of \eqref{-211}. The momentum equation \eqref{-217} and the induction equation \eqref{-218} hold true according to \eqref{-224} and \eqref{-186}. The initial conditions \eqref{-219} follow as in the previous limit passages, in particular the initial conditions for $\rho u$ and $B$ can be derived as the one for the magnetic induction on the $\eta$-level in Section \ref{inductionequationetalevel}. In the energy inequality \eqref{-167} on the $\alpha$-level we can pass to the limit using the weak lower semicontinuity of norms to infer the energy inequality \eqref{-220}. Finally, the compatibility of $u$ with $\lbrace S_0^i, X^i \rbrace_{i=1}^m$ is shown in \eqref{-314}. This finishes the proof of Theorem \ref{mainresult}.

{\centering \section{Appendix} \par }

For the construction and estimation of the density on the $\Delta t$-level we use some classical results for the regularized continuity equation which are summarized in the following lemma.

\begin{satz}
\label{neumannproblem}
Let $T>0$, $r \in (0,1)$ and assume $\Omega \subset \mathbb{R}^3$ to be a bounded domain of class $C^{2,r}$. Let $V_n$ denote the Galerkin space from our approximate system in Section \ref{approximatesystem} and let $w \in C([0,T];V_n)$. Finally, consider the initial data $\rho_0 \in C^{2,r}(\overline{\Omega})$ such that
\begin{align}
\left. \nabla \rho_0 \cdot \text{n} \right|_{\partial \Omega} = 0,\quad \quad \underline{\rho} \leq \rho_0 \leq \overline{\rho}\quad \text{in } \overline{\Omega} \nonumber
\end{align}

for two constants $0 < \underline{\rho} \leq \overline{\rho} < \infty$. Then the Neumann problem
\begin{align}
\partial_t \rho + \operatorname{div}(\rho w) = \epsilon \Delta \rho \quad \text{in } (0,T)\times \Omega, \quad \quad \left. \nabla \rho_0 \cdot \operatorname{n} \right|_{\partial \Omega} = 0,\quad \quad \rho(0) = \rho_0 \quad \text{in } \Omega \label{-328}
\end{align}

admits a unique solution $\rho = \rho(w)$ in the class
\begin{align}
\rho \in C \left([0,T];C^{2,r}\left(\overline{\Omega}\right)\right) \bigcap C^1 \left([0,T];C^{0,r}\left(\overline{\Omega}\right)\right). \label{-332}
\end{align}

In addition, the estimates
\begin{align}
0 < \underline{\rho} \exp \left(-\int_{0}^t \left\| w(\tau) \right\|_{W^{1,\infty}(\Omega)}\ d\tau \right) \leq \rho(w)(t) \leq \overline{\rho} \exp \left(\int_{0}^t \left\| w(\tau) \right\|_{W^{1,\infty}(\Omega)}\ d\tau \right) < \infty \quad &\forall t \in [0,T], \label{-329} \\
\left\| \rho(w) \right\|_{C \left([0,T];C^{2,r}\left(\overline{\Omega}\right)\right)} + \left\| \rho(w) \right\|_{C^1 \left([0,T];C^{0,r}\left(\overline{\Omega}\right)\right)} \leq c(w)\quad \forall w \in C([0,T];V_n) & \label{-330} \\
\left\| \rho \left(w_1\right) - \rho \left(w_2\right) \right\|_{C([0,T];L^2(\Omega))} \leq c \left\| w_1 - w_2 \right\|_{C([0,T];W^{1,\infty}(\Omega))} \quad \forall w_1,w_2 \in C([0,T];&V_n) \label{-331}
\end{align}

are satisfied for some constant $c(w)>0$ bounded on bounded subsets of $C([0,T];V_n)$ and some constant $c>0$ independent of $w_1$, $w_2$.
\end{satz}

\textbf{Proof}

The existence of a unique solution $\rho$ to the Neumann problem \eqref{-328} in the class
\begin{align}
\rho \in C\left([0,T];W^{1,2}(\Omega) \right) \bigcap L^2\left(0,T;W^{2,2}(\Omega) \right),\quad \partial_t \rho \in L^2\left((0,T)\times \Omega \right), \nonumber
\end{align}

which satisfies the estimates \eqref{-329} and \eqref{-331}, is well known, c.f.\ \cite[Proposition 7.39]{novotnystraskraba}. The additional regularity \eqref{-332} together with the corresponding estimate \eqref{-330} then follows by classical results (\cite[Theorem 10.22, Theorem 10.23]{singularlimits}) on the maximal regularity for parabolic problems, c.f. \cite[Lemma 3.1]{singularlimits}.

$\hfill \Box$

For the limit passage in the time discretization we use a modified version of \cite[Theorem 8.9]{roubicek} to infer that different interpolants of the discrete quantities converge to the same limit function:

\begin{satz}
\label{equalityofrothelimits}
Let $\Omega \subset \mathbb{R}^3$ be a bounded domain and assume that
\begin{align}
f_{\Delta t} \buildrel\ast\over\rightharpoonup f \quad \text{in } L^\infty(0,T;L^2(\Omega)),\quad \overline{f}_{\Delta t} \buildrel\ast\over\rightharpoonup \overline{f} \quad \text{in } L^\infty(0,T;L^2(\Omega))\quad \overline{f}'_{\Delta t} \buildrel\ast\over\rightharpoonup \overline{f}' \quad \text{in } L^\infty(0,T;L^2(\Omega)), \nonumber
\end{align}

where $f_{\Delta t}, \overline{f}_{\Delta t}, \overline{f}_{\Delta t}'$ denote piecewise affine and piecewise constant interpolants of some discrete functions $f^k_{\Delta t}$, $k=0,...,\frac{T}{\Delta t}$ as introduced in (\ref{21})--(\ref{23}). Then the limit functions $f,\ \overline{f},\ \overline{f}'$ coincide,

\begin{align}
f = \overline{f} = \overline{f}'\quad \text{a.e. in } (0,T)\times \Omega. \label{918}
\end{align}

\end{satz}

\textbf{Proof}

The proof of the first identity in \eqref{918} is given in \cite[Theorem 8.9]{roubicek}, a simple modification then also yields the second equality, c.f.\ \cite[Lemma 7.1]{incompressiblecase}.

$\hfill \Box$

In order to show convergence of the positions of the solid bodies, we make use of the following lemma which can be found in \cite[Proposition 5.1]{feireisl}.

\begin{satz} \label{convergence}

Let $O \subset \mathbb{R}^3$ be a bounded domain of class $C^2\bigcap C^{0,1}$ and let $\delta > 0$. Let further $\lbrace u_k \rbrace_{k \in \mathbb{N}}$ be a sequence of vector fields bounded in $L^2(0,T;W^{1,\infty}(\mathbb{R}^3))$ uniformly with respect to $k$. Moreover, denote by $X_k$ the Carathéodory solution to the initial value problem
\begin{align}
\frac{dX_k(t;x)}{dt} =  u_k \left( t, X_k(t;x) \right),\quad \text{for } t \in [0,T], \quad \quad X_k(0;x) =  x,\quad \text{for } x \in \mathbb{R}^3 \label{-338}
\end{align}

and set $O_k(t):= X_k(t;O)$ and $S_k(t):= (O_k(t))^\delta$. Then there exists a function $X: [0,T]\times \mathbb{R}^3 \rightarrow \mathbb{R}^3$ such that, for a suitable subsequence, it holds
\begin{align}
X_k \rightarrow& X \quad &\text{in }& C\left([0,T];C_{\operatorname{loc}} \left( \mathbb{R}^3 \right) \right),\label{-334} \\
\textbf{db}_{O_k(t)} \rightarrow& \textbf{db}_{O(t)} \quad &\text{in }& C_{\operatorname{loc}} \left( \mathbb{R}^3 \right), \label{-335} \\
\textbf{db}_{S_k(t)} \rightarrow& \textbf{db}_{S(t)} \quad &\text{in }& C_{\operatorname{loc}} \left( \mathbb{R}^3 \right) \label{-336}
\end{align}

uniformly with respect to $t \in [0,T]$, where $O(t):= X(t;O)$ and $S(t):= (O(t))^\delta$. Moreover, if the extracted subsequence is chosen such that
\begin{align}
u_k \buildrel\ast\over\rightharpoonup u \quad \text{in } L^2\left(0,T;W^{1,\infty}\left(\mathbb{R}^3\right)\right) \nonumber
\end{align}

for some $u \in L^2(0,T;W^{1,\infty}(\mathbb{R}^3))$, then $X$ constitutes the Carathéodory solution to
\begin{align}
\frac{dX(t;x)}{dt} =  u \left( t, X(t;x) \right),\quad \text{for } t \in [0,T], \quad \quad X(0;x) =  x,\quad \text{for } x \in \mathbb{R}^3. \label{-337}
\end{align}

\end{satz}

\textbf{Proof}

The convergences \eqref{-334} and \eqref{-335} can be concluded via an application of the Arzelà-Ascoli theorem, the convergence \eqref{-336} follows directly from the convergence $\eqref{-335}$. The relation \eqref{-337} then follows by passing to the limit in the corresponding relation \eqref{-338}. For the details of the proof we refer to \cite[Proposition 5.1]{feireisl}.

$\hfill \Box$

Since, according to Lemma \ref{compatibility}, the rigid bodies can never leave the domain $\Omega$, we moreover have the following corollary for the case that the velocity fields $u_k$ in Lemma \ref{convergence} are compatible with suitable systems of isometries.

\begin{corollary} \label{convergenceandcompatibility}

Let $\Omega, S_0^1,...,S_0^m \subset \mathbb{R}^3$, $m \in \mathbb{N}$, be bounded domains of class $C^2\bigcap C^{0,1}$. Let further $\lbrace u_k \rbrace_{k \in \mathbb{N}}$ be a sequence of vector fields bounded in $L^2(0,T;H^{1,2}(\Omega))$ uniformly with respect to $k$ and let each $u_k$ be compatible with the system $\lbrace S_0^i,X_k^i \rbrace_{i=1}^m$ where $X_k^i(t):\mathbb{R}^3 \rightarrow \mathbb{R}^3$, $t \in [0,T]$, $i=1,...,m$ denotes an isometry. Then there exist isometries $X^i(t):\mathbb{R}^3 \rightarrow \mathbb{R}^3$ such that, for a suitable subsequence, it holds
\begin{align}
X_k^i \rightarrow X^i \quad \text{in } C\left([0,T];C_{\operatorname{loc}} \left( \mathbb{R}^3 \right) \right) \nonumber
\end{align}

and if the extracted subsequence is chosen such that
\begin{align}
u_k \buildrel\ast\over\rightharpoonup u \quad \text{in } L^2\left(0,T;H^{1,2}(\Omega)\right) \nonumber
\end{align}

for some $u \in L^2(0,T;H^{1,2}(\Omega))$, then $u$ is compatible with the system $\lbrace S_0^i,X_k^i \rbrace_{i=1}^m$.

\end{corollary}

\begin{center}
\Large\textbf{Acknowledgements} \\[4mm]
\end{center}

This article is part of the PhD project of the author. J.S. would like to thank Barbora Benešová, Šárka Nečasová and Anja Schlömerkemper for many fruitful discussions. This work has been supported by the Czech Science Foundation (GA\v CR) through the project 22-08633J. The work was moreover supported by the Charles University, Prague via project PRIMUS/19/SCI/01.

\end{document}